\documentclass[11pt,colorlinks,urlcolor=cyan,citecolor=MidnightBlue,linkcolor=MidnightBlue,dvipsnames]{article}
\usepackage{t1enc}
\usepackage{lmodern}
\usepackage[T2A,T1]{fontenc}
\usepackage[utf8]{inputenc}
\usepackage[russian,english]{babel}
\usepackage{amsmath,amssymb,amsfonts,amsthm,mathrsfs,textcomp,url,bbm,cancel,comment,enumerate,orcidlink}
\usepackage{graphicx}
\usepackage{float}
\usepackage[all,ps]{xy}
\usepackage[left=3cm,right=3cm,top=2.8cm,bottom=2.8cm]{geometry}
\usepackage{longtable}
\usepackage{tikz}
\usepackage{nicematrix}
\makeatletter
\newcommand{\affil}[3]{\gdef\@affila{\textsc{#1}}
  \gdef\@affilb{\textsc{#2}}
  \gdef\@affilc{\textsc{#3}}}
\newcommand{\address}[3]{\gdef\@addressa{#1}
  \gdef\@addressb{#2}
  \gdef\@addressc{#3}}
\newcommand{\email}[3]{\gdef\@emaila{\url{#1}}
  \gdef\@emailb{\url{#2}}
  \gdef\@emailc{\url{#3}}}
\newcommand{\@endstuff}{\par\vspace{\baselineskip}\noindent\footnotesize
  \begin{tabular}{@{}l}
    \@affila\\
    \@addressa\\
    \textit{E-mail address:} \@emaila
  \end{tabular}
  \par\medskip\noindent
  \begin{tabular}{@{}l}
    \@affilb\\
    \@addressb\\
    \textit{E-mail address:} \@emailb
  \end{tabular}
  \par\medskip\noindent
  \begin{tabular}{@{}l}
    \@affilc\\
    \@addressc\\
    \textit{E-mail address:} \@emailc
  \end{tabular}}
\AtEndDocument{\@endstuff}
\makeatother

\usetikzlibrary{matrix,calc,decorations,positioning}
\pgfkeys{/tikz/.cd,
  alt double distance/.initial=5pt,
  alt double step/.initial=1pt,}
\pgfdeclaredecoration{double deco}{initial}
{
  \state{initial}[width=\pgfkeysvalueof{/tikz/alt double step},next state=cont] {
    \pgfmoveto{\pgfpoint{\pgfkeysvalueof{/tikz/alt double step}}{\pgfkeysvalueof{/tikz/alt double distance}/2}}
    \pgfpathlineto{\pgfpoint{0.3\pgflinewidth}{\pgfkeysvalueof{/tikz/alt double distance}/2}}
    \pgfpathmoveto{\pgfpoint{0.3\pgflinewidth}{-\pgfkeysvalueof{/tikz/alt double distance}/2}}
    \pgfpathlineto{\pgfpoint{1pt}{-\pgfkeysvalueof{/tikz/alt double distance}/2}}
    \pgfcoordinate{lastup}{\pgfpoint{1pt}{\pgfkeysvalueof{/tikz/alt double distance}/2}}
    \pgfcoordinate{lastdown}{\pgfpoint{1pt}{-\pgfkeysvalueof{/tikz/alt double distance}/2}}}
  \state{cont}[width=\pgfkeysvalueof{/tikz/alt double step}]{
    \pgfmoveto{\pgfpointanchor{lastup}{center}}
    \pgfpathlineto{\pgfpoint{\pgfkeysvalueof{/tikz/alt double step}}{\pgfkeysvalueof{/tikz/alt double distance}/2}}
    \pgfcoordinate{lastup}{\pgfpoint{\pgfkeysvalueof{/tikz/alt double step}}{\pgfkeysvalueof{/tikz/alt double distance}/2}}
    \pgfmoveto{\pgfpointanchor{lastdown}{center}}
    \pgfpathlineto{\pgfpoint{\pgfkeysvalueof{/tikz/alt double step}}{-\pgfkeysvalueof{/tikz/alt double distance}/2}}
    \pgfcoordinate{lastdown}{\pgfpoint{\pgfkeysvalueof{/tikz/alt double step}}{-\pgfkeysvalueof{/tikz/alt double distance}/2}}}
  \state{final}[width=0pt]
  { 
    \pgfmoveto{\pgfpointdecoratedpathlast}}}
\tikzset{alt double/.style={decorate,decoration=double deco}}
\allowdisplaybreaks

\newcommand{\toverset}[2]{%
  \mathop{#2}\limits^{\vbox to -.1ex{\kern-0.4ex\hbox{$\scriptstyle #1$}\vss}}}
\newcommand{\tightoverset}[2]{%
  \mathop{#2}\limits^{\vbox to -.5ex{\kern-0.4ex\hbox{$\scriptstyle #1$}\vss}}}
\renewcommand{\underset}[2]{%
  \mathop{#2}\limits_{\vbox to -.5ex{\kern-1.6ex\hbox{$\scriptstyle #1$}\vss}}}
\newcommand{\tightunderset}[2]{%
  \mathop{#2}\limits_{\vbox to -.5ex{\kern-1.8ex\hbox{$\scriptstyle #1$}\vss}}}

\def\Z{\mathbb{Z}}

\def\R{\mathbb{R}}

\def\into{\hookrightarrow}
\def\imto{\looparrowright}

\def\ol{\overline}

\def\la{\langle}
\def\ra{\rangle}

\def\im{\mathop{\rm im}}
\def\coker{\mathop{\rm coker}}

\def\id{\mathop{\rm id}}
\def\pr{\mathop{\rm pr}}

\def\Tp{\mathrm{Tp}}

\def\Sq{\mathrm{Sq}}

\def\O{\mathrm{O}}

\def\SO{\mathrm{SO}}

\def\Hom{\textstyle{\mathop{\rm Hom}}}

\newenvironment{thm'}[1]%
{\phantomsection\par\addvspace{.5\baselineskip}\noindent\textbf{Theorem \ref*{#1}'\label{#1'}.\enspace\ignorespaces}\begin{em}}%
  {\end{em}\par\addvspace{.5\baselineskip}}%
\newenvironment{prf}[1][\unskip]%
{\par\addvspace{.5\baselineskip}\noindent\textbf{Proof #1.\enspace\ignorespaces}}%
{~$\square$\par\addvspace{.5\baselineskip}}%
{\par\addvspace{.5\baselineskip}\noindent\textbf{Sketch of the proof #1.\enspace\ignorespaces}}%
{~$\square$\par\addvspace{.5\baselineskip}}%
{\par\addvspace{.5\baselineskip}\noindent\textit{Claim.\enspace\ignorespaces}\begin{em}}%
  {\end{em}\par\addvspace{.5\baselineskip}}%
{\par\addvspace{.5\baselineskip}\noindent\textit{Proof #1.\enspace\ignorespaces}}%
{~$\diamond$\par\addvspace{.5\baselineskip}}%
{\par\addvspace{.5\baselineskip}\noindent\textit{Remark.\enspace\ignorespaces}}%
{\par\addvspace{.5\baselineskip}}%
\newenvironment{subjclass}[1]%
{\par\noindent\begin{small}\textit{#1 Mathematics Subject Classification.\enspace\ignorespaces}}%
  {\end{small}\par}%
\newenvironment{key}%
{\par\noindent\begin{small}\textit{Key words and phrases.\enspace\ignorespaces}}%
  {\end{small}\par}%
{\par\noindent\begin{small}\textbf{Acknowledgement.\enspace\ignorespaces}}%
  {\end{small}\par}%
{\par\addvspace{.5\baselineskip}\noindent\begin{em}}%
  {\end{em}\par\addvspace{.5\baselineskip}}%
\newtheorem{thm}{Theorem}%
\newtheorem{lemma}{Lemma}%
\newtheorem{crly}{Corollary}%
\newtheorem*{ques}{Question}%
\newtheorem*{conj}{Conjecture}%
\theoremstyle{definition}
\newtheorem{defi}{Definition}%
\newtheorem*{rmk}{Remark}%
\newtheorem*{ex}{Example}%
\newcommand{\ZeroRoman}[1]{%
\ifcase\value{#1}\relax 0\else\Roman{#1}\fi}

\newcounter{t}

\numberwithin{equation}{section}
\title{On coincidences of Morin and first order Thom--Boardman singular loci}
\author{András Csépai \and András Szűcs \and Tamás Terpai \orcidlink{0000-0002-5707-2668}}
\affil{ELTE Eötvös Loránd University, Budapest, Hungary, Institute of Mathematics}{ELTE Eötvös Loránd University, Budapest, Hungary, Institute of Mathematics}{ELTE Eötvös Loránd University, Budapest, Hungary, Institute of Mathematics}
\address{Pázmány Péter sétány 1/c, Budapest, H-1117 Hungary}{Pázmány Péter sétány 1/c, Budapest, H-1117 Hungary}{Pázmány Péter sétány 1/c, Budapest, H-1117 Hungary}
\email{csepai.andras112358@gmail.com}{andras.szucs@ttk.elte.hu}{terpai@math.elte.hu}
\date{}
\begin{document}

\maketitle

\begin{abstract}
  It is well-known that the Thom polynomial in Stiefel--Whitney classes expressing the cohomology class dual to the locus of the cusp singularity for codimension-$k$ maps and that of the corank-$2$ singularity for codimension-$(k-1)$ maps coincide. The aim of the present paper is to find out whether there is any geometric explanation to this seemingly mysterious coincidence. We thank László Fehér for posing us this interesting question that we answer here in the positive, and motivated by it we search for further similar coincidences of the loci of the corank-$r$ and Morin singular points, but found such only for special classes of maps. Finally we compute the cohomology classes dual to the singularity strata for any Morin map. (This result--for holomorphic maps--was presented by Kazarian but the proof was left unpublished.)

\end{abstract}

\begin{subjclass}{2020}
  57R45 (Primary); 57R70; 57R20 (Secondary)
\end{subjclass}

\begin{key}
  Thom polynomials; relation of singular loci
\end{key}

\section{Introduction}
\label{sec:intro}

The present paper investigates the relation and possible coincidence up to homology or cobordism (or other equivalences) of certain singular loci for smooth (real $C^\infty$) maps, and the geometric reasons behind such relations. More concretely, we will study the loci of the singularities $\Sigma^{1_r}(k)$ and $\Sigma^r(l)$ where $\Sigma^{1_r}$ is short for $\Sigma^{1,\ldots,1,0}$ with $r$ digits $1$ and the arguments in parentheses denote the codimensions of the maps under consideration. The main motivation for this is the following:

\begin{ques}
  Classical computations show that the (real, modulo $2$) Thom polynomials of the singularities $\Sigma^{1_2}(k)$ and $\Sigma^2(k-1)$ coincide, they are both the polynomial $w_{k+1}^2+w_kw_{k+2}$. Is there a geometric explanation for this 
  coincidence?
\end{ques}

This question was posed 
by L. Fehér and one purpose of this paper is 
to answer it to the affirmative (see theorem \ref{thm:cusp} and corollary \ref{crly:cusp}). It also motivated us to investigate the relation of the loci of the singularities $\Sigma^{1_r}(k)$ and $\Sigma^r(k-r+1)$ for any $r$ for special types of maps (see theorems \ref{thm:prim}, \ref{thm:twprim} and corollaries \ref{crly:prim}, \ref{crly:twprim}). This is related to the problem of computing the Thom polynomials of the singularities $\Sigma^{1_r}(k)$ for general Morin maps (see theorem \ref{thm:mor}) which was also studied by Kazarian \cite{morchar}. Before stating our results mentioned above, 
let us briefly introduce the objects that appear in their formulations. 


A central problem in global singularity theory is, for a given singularity type $\eta$ of smooth map germs, to understand the singular locus $\eta(f)\subset M$ of maps $f\colon M\to N$, that is, the points where the singularity of $f$ is $\eta$. A classical result of Thom \cite{tp} 
is that 
there
is an element $\Tp_{\eta(k)}\in\Z_2[w_1,w_2,\ldots]\cong H^*(B\O;\Z_2)$, called now the \textit{Thom polynomial}, such that for any generic smooth map $f\colon M\to N$ of codimension $k$, the Poincaré dual of the closure $\ol\eta(f)$ of the $\eta$-locus of $f$ 
is obtained by substituting the Stiefel--Whitney classes of the stable normal bundle $\nu_f:=f^*TN\ominus TM$ in the variables of $\Tp_{\eta(k)}$. 
In some cases (characterised by Vassiliev \cite{vas}) the Thom polynomial of $\eta(k)$ also has an oriented analogue $\Tp^\SO_{\eta(k)}$ in the ring $H^*(B\SO;\Z)$ giving the dual of $\ol\eta(f)$ in integral cohomology. Note that $H^*(B\SO;\Z)$ is of the form $\Z[p_1,p_2,\ldots]\oplus\im\beta$ where the summand $\im\beta$ is isomorphic to the image of the Steenrod operation $\Sq^1$; for any element $w_{i_1}\ldots w_{i_r}\in\im\Sq^1$ we will denote by $v_{i_1,\ldots,i_r}$ the corresponding element in $\im\beta$.

Our goal will be to understand geometrically the coincidences of Thom polynomials 
for maps of different codimensions. This presents a problem which we overcome by introducing the following notion.

\begin{defi}
  For a smooth map $f\colon M\to N$, an \textit{($l$-tuple) de-suspension} of $f$ is a generic perturbation $\tilde f\colon M\times\R^l\to N$ of the composition $f\circ\pr_M$ (i.e. a generic map close to it in the space of smooth maps) where $\pr_M\colon M\times\R^l\to M$ is the projection to $M$, such that the submanifold $M\times\{0\}$ is transverse to the singular strata of $\tilde f$. (Note that the de-suspension of a map is not unique, but it is unique homotopically.)
\end{defi}

Our main results will geometrically relate $\Sigma^{1_r}$-loci of maps with $\Sigma^r$-loci of their de-suspensions for three special types of maps defined as follows.

\begin{defi}
  A generic map $f\colon M\to N$ is said to be a \textit{Morin map} if $\dim\ker df_p$ (the corank of its differential) is at most $1$ for all $p\in M$. The Morin map $f$ will be called
  \begin{enumerate}
  \item a \textit{cusp map} if it only has singularities of types $\Sigma^{1_r}$ with $r\le2$, i.e. all of its singular points belong to the \textit{fold} ($\Sigma^{1,0}$) or the \textit{cusp} ($\Sigma^{1,1,0}$) locus, 
  \item a \textit{prim map} if the kernel line bundle $\ker df$ over the set $\Sigma(f)$ of singular points is trivial,
  \item a \textit{twisted prim map} if there is a line bundle $\ell$ over $M$ such that the kernel bundle $\ker df$ is the restriction $\ell|_{\Sigma(f)}$.
  \end{enumerate}
\end{defi}

\begin{rmk}
  The name ``prim'' is a shortening of \textit{pr}ojected \textit{im}mersion. It is not hard to see that prim maps are precisely the hyperplane projections of immersions, that is, a generic map $f\colon M\to N$ is prim precisely if there is an immersion $j\colon M\imto N\times\R$ such that $f$ is the composition $\pr_N\circ j$; see e.g. \cite{immemb}. 
  Such a $j$ will be called an \textit{immersion lift} of $f$.
\end{rmk}

The paper is organised as follows: in section \ref{sec:result} we state our main results; in section \ref{sec:pre} we define the notions of singularities and Thom polynomials used throughout the paper and recall their basic properties; in section \ref{sec:prf} we prove theorems \ref{thm:cusp}, \ref{thm:prim}, \ref{thm:twprim}, \ref{thm:mor} and corollaries \ref{crly:cusp}, \ref{crly:prim}, \ref{crly:twprim} below; and finally in section \ref{sec:fin} we pose open problems related to the topic. We also included an \hyperref[sec:app]{appendix} containing tedious but crucial computations needed for the proof of theorem \ref{thm:cusp}.

\section{Formulation of the results}
\label{sec:result}

Relating the loci of the singularity types $\Sigma^{1_2}(k)$ and $\Sigma^2(k-1)$ we have:

\begin{thm}
  \label{thm:cusp}
  If $f\colon M^n\to N^{n+k}$ is a cusp map where $M$ is closed and $k\ge1$ and $\tilde f\colon M\times\R\to N$ is a de-suspension of $f$, then the submanifolds $\Sigma^{1_2}(f)$ and $\Sigma^2(\tilde f)\cap M\times\{0\}$ of $M$ are cobordant as embedded submanifolds in $M$. This cobordism is oriented if $M$ and $N$ are oriented and the codimension $k$ is odd.
\end{thm}

We note here that the perturbation $\tilde f$ of the map $f\circ\pr_M$ that we construct in the proof of theorem \ref{thm:cusp} will be such that the submanifolds $\Sigma^{1_2}(f)$ and $\Sigma^2(\tilde f)\cap M\times\{0\}$ actually coincide. From this theorem we deduce:


\begin{crly}
  \label{crly:cusp}
  We have
  $$\Tp_{\Sigma^{1_2}(k)}=\Tp_{\Sigma^2(k-1)}$$
  and if $k$ is odd, then we also have
  $$\Tp^\SO_{\Sigma^{1_2}(k)}=\Tp^\SO_{\Sigma^2(k-1)}.$$
\end{crly}

\begin{rmk}
  The first equation in corollary \ref{crly:cusp} is, as noted before, well-known. The second equation is apparently not so widely known, but it also simply follows from known results: using the formulae in \cite{calcsk} and \cite{singprim} (see also \cite{tps11}) we get $\Tp^\SO_{\Sigma^{1_2}(k)}=p_{\frac{k+1}2}+v_{k,k+2}=\Tp^\SO_{\Sigma^2(k-1)}$ for $k$ odd. Note that if $k$ is even, then these singular loci do not represent integral cohomology classes in general, hence then the Thom polynomials $\Tp^\SO_{\Sigma^{1_2}(k)}$ and $\Tp^\SO_{\Sigma^2(k-1)}$ are not defined.
\end{rmk}

As corollary \ref{crly:cusp} shows, theorem \ref{thm:cusp} gives a geometric reason for the coincidence of the Thom polynomials of the singularities $\Sigma^{1_2}(k)$ and $\Sigma^2(k-1)$ 
through the embedded cobordism of singular loci. We will now consider two special classes of maps for which analogous relations hold for the loci of the singularities $\Sigma^{1_r}(k)$ and $\Sigma^r(k-r+1)$ for any $r$. We also note in advance that these singular loci represent integral cohomology classes for general maps only when $k$ is odd and $r$ is even.

\begin{thm}
  \label{thm:prim}
  If $f\colon M^n\to N^{n+k}$ is a prim map where $M$ is closed and $k\ge r-1$ and $\tilde f\colon M\times\R^{r-1}\to N$ is a de-suspension of $f$, then the submanifolds $\ol\Sigma^{1_r}(f)$ and $\Sigma^r(\tilde f)\cap M\times\{0\}$ of $M$ are cobordant as embedded submanifolds in $M$. This cobordism is oriented if $M$ and $N$ are oriented and the codimension $k$ is odd.
\end{thm}

This yields equations of Thom polynomials similar to the ones in corollary \ref{crly:cusp} evaluated for $k$-codimensional prim maps.

\begin{crly}
  \label{crly:prim}
  If $f\colon M^n\to N^{n+k}$ is a prim map where $M$ is closed, then for any $r$ we have
  $$\Tp_{\Sigma^{1_r}(k)}(f)=\Tp_{\Sigma^r(k-r+1)}(f)=w_{k+1}(\nu_f)^r$$
  and if $M$ and $N$ are oriented and $k$ is odd and $r$ is even, then we also have
  $$\Tp^\SO_{\Sigma^{1_r}(k)}(f)=\Tp^\SO_{\Sigma^r(k-r+1)}(f)=p_{\frac{k+1}2}(\nu_f)^{\frac r2}.$$
\end{crly}

To explain the second term in both equations above, note that although a $k$-codimensional prim map cannot have $\Sigma^r(k-r+1)$-singularities, it still makes sense to substitute the characteristic classes of its stable normal bundle into the variables of the abstract Thom polynomials of $\Sigma^r(k-r+1)$.

\begin{rmk}
  The first formula in corollary \ref{crly:prim} is in essence not new: the classical Giambelli--Thom--Porteous formula implies that the leading term of $\Tp_{\Sigma^r(k-r+1)}$ is $w_{k+1}^r$ and all other terms contain a factor $w_l$ with $l>k+1$; the same is proved for $\Tp_{\Sigma^{1_r}(k)}$ in \cite{rim}\footnote{Most papers on Thom polynomials consider them in the holomorphic setting but a classical result of Borel and Haefliger \cite{borhae} implies that the smooth setting (with $\Z_2$ coefficients) is straightforward from that.} and since for a codimension-$k$ prim map $f$ the Stiefel--Whitney classes $w_l(\nu_f)$ vanish for $l>k+1$, the statement follows.
\end{rmk}

Theorem \ref{thm:prim} and corollary \ref{crly:prim} can, in a weaker form, be extended from prim maps to the more general class of twisted prim maps.

\begin{thm}
  \label{thm:twprim}
  If $f\colon M^n\to N^{n+k}$ is a twisted prim map where $M$ is closed and $k\ge r-1$ and $\tilde f\colon M\times\R^{r-1}\to N$ is a de-suspension of $f$, further, $M$ and $N$ are oriented, $k$ is odd and $r$ is even, then the homology classes $2[\ol\Sigma^{1_r}(f)]$ and $2[\Sigma^r(\tilde f)\cap M\times\{0\}]$ in $H_*(M;\Z)$ coincide.
\end{thm}

\begin{crly}
  \label{crly:twprim}
  If $f\colon M^n\to N^{n+k}$ is a twisted prim map where $M$ is closed, $M$ and $N$ are oriented and $k$ is odd and $r$ is even, then we have
  $$2\Tp^\SO_{\Sigma^{1_r}(k)}(f)=2\Tp^\SO_{\Sigma^r(k-r+1)}(f).$$
\end{crly}

Finally, related to the connection of the Thom polynomials of the singularities $\Sigma^{1_r}$ and $\Sigma^r$, we give a formula computing the Thom polynomials of all singularities $\Sigma^{1_r}$ restricted to Morin maps.

\begin{rmk}
  The Thom polynomial $\Tp_{\Sigma^{1_r}(k)}$ is not known in general, it has been computed for a general $k$ only if $r$ is at most $6$; see \cite{tpmor}. In contrast, the Thom polynomial $\Tp_{\Sigma^r(k)}$ is classically known for all $k$ and $r$, it is determined by the Giambelli--Thom--Porteous formula; see e.g. \cite{gtp}.
\end{rmk}

\begin{thm}
  \label{thm:mor}
  If $f\colon M^n\to N^{n+k}$ is a Morin map, 
  then we have
  $$\Tp_{\Sigma^{1_r}(k)}(f)=
  \begin{cases}
    \left(w_{k+1}(\nu_f)^2+w_k(\nu_f)w_{k+2}(\nu_f)\right)^{\frac r2},&\text{if }r\text{ is even}\\
    w_{k+1}(\nu_f)\left(w_{k+1}(\nu_f)^2+w_k(\nu_f)w_{k+2}(\nu_f)\right)^{\frac{r-1}2},&\text{if }r\text{ is odd}
  \end{cases}
  $$
  and if $M$ and $N$ are oriented and $k$ is odd and $r$ is even, then we also have
  $$\Tp^\SO_{\Sigma^{1_r}(k)}(f)=\left(p_{\frac{k+1}2}(\nu_f)+v_{k,k+2}(\nu_f)\right)^{\frac r2}.$$
\end{thm}

\begin{rmk}
  The first formula in theorem \ref{thm:mor} follows from a result of Kazarian \cite{morchar} where this Thom polynomial was computed for holomorphic Morin maps. Our proof of it is the same as Kazarian's, but we will still repeat it since it fits with the techniques of the present paper and although Kazarian's formula was published in \cite{kaztp}, to our knowledge its proof remained unpublished.
\end{rmk}

\section{Singularities and Thom polynomials}
\label{sec:pre}


Throughout the paper we consider singular loci of generic smooth maps of non-negative codimension. The notions of singularity and genericity are defined in various ways depending on context, in the present paper ``singularity'' will always mean Thom--Boardman singularity 
which we define in the following (see also \cite{singdiff}) together with our definition of the notion of genericity.

\begin{defi}
  \label{defi:tbsing}
  For any map $f\colon M^n\to N^{n+k}$ and number $r$ we define $\Sigma^r(f)\subset M$ to be the set of points $p\in M$ where the corank of the differential of $f$ (i.e. $\dim\ker df_p$) is $r$; the set of \textit{singular points} of $f$ is $\Sigma(f):=\Sigma^1(f)\cup\ldots\cup\Sigma^n(f)$. The subset $\Sigma^r(f)$ is the preimage under the $1$-jet extension $j^1f\colon M\to J^1(M,N)$ of the subset $\Sigma^r(M,N)\subset J^1(M,N)$ of $1$-jets of corank $r$ which is a submanifold of this jet bundle. Now if $j^1f$ is transverse to the submanifold $\Sigma^r(M,N)$, then the subset $\Sigma^r(f)\subset M$ is a manifold, hence we can consider for any $s\le r$ the subset $\Sigma^{r,s}(f):=\Sigma^s(f|_{\Sigma^r(f)})$ which is the preimage under the $2$-jet extension $j^2f$ of the correspondingly defined submanifold $\Sigma^{r,s}(M,N)\subset J^2(M,N)$. Iterating this process yields for any weakly decreasing sequence $S$ of non-negative integers the submanifold $\Sigma^S(f)\subset M$ called a \textit{Thom--Boardman stratum of $f$}; this stratum is again the $j^{|S|}f$-preimage of the \textit{Thom--Boardman stratum} $\Sigma^S(M,N)\subset J^{|S|}(M,N)$. We call the map $f$ \textit{generic} if its jet extensions are transverse to all Thom--Boardman strata in the corresponding jet bundles. The points in $\Sigma^S(f)$ are called the \textit{$\Sigma^S$-points} of $f$ and $\Sigma^S(f)$ is also called the \textit{$\Sigma^S$-locus} of $f$; this classification of singular points for all maps $f\colon M^n\to N^{n+k}$ where $n$ is arbitrary and $k$ is fixed defines the \textit{Thom--Boardman singularity type} $\Sigma^S(k)$.
\end{defi}


\begin{rmk}
  As seen above, we denote by $\eta(k)$ a singularity type $\eta$ considered for maps of fixed codimension $k$, on the other hand we denote by $\eta(f)$ the $\eta$-locus of a map $f$. Although this makes the notations ambiguous, it should not be confusing.
\end{rmk}


\begin{rmk}
  If $S$ is a decreasing sequence, then for a map $f\colon M\to N$ and a point $p\in M$ the local germ of $f$ at $p$ determines whether $p$ is a $\Sigma^S$-point or not, but the germs of $f$ at different points of $\Sigma^S(f)$ may be different up to local coordinate change. However, Morin's classical result \cite{mor} shows that if $f\colon M^n\to N^{n+k}$ is a Morin map (i.e. if $\Sigma(f)= \Sigma^1(f)$), then the germs of $f$ at two points of $M$ are the same in some local coordinate systems precisely if they belong to the same Thom--Boardman singularity $\Sigma^{1_r}(k)$; this singularity is often denoted by $A_r(k)$ and also called a \textit{Morin singularity}. Thus the Morin singularities are locally well understood, however, this is not the case for their global structure.
\end{rmk}


As noted in the introduction, for a singularity $\eta(k)$ the Thom polynomial $\Tp_{\eta(k)}$ describes the closure of the $\eta$-locus of any generic $k$-codimensional map $f\colon M\to N$ up to (co)homology. We note here that generally the closure $\ol\eta(f)\subset M$ is not a submanifold, however, it is in two special cases that are important in this paper: firstly if $f$ is a Morin map, then the closure of $\Sigma^{1_r}(f)=\Sigma^{1,\ldots,1,0}(f)$ (with $r$ digits $1$) is $\Sigma^{1,\ldots,1}(f)$ (again with $r$ digits $1$) which is a submanifold; secondly if $\Sigma^{r+1}(f)$ is empty, then $\Sigma^r(f)$ is a closed submanifold.

Modern approaches for constructing the Thom polynomial $\Tp_{\eta(k)}$ are due to Kazarian \cite{kazspace} and Rimányi \cite{rim}, we now briefly recall the construction of it based on \cite{kazspace}. Consider the infinite jet bundle $J^\infty_0(\varepsilon^l,\gamma_{l+k})$ over $B\O(l+k)$ where $l$ is any sufficiently large integer, $\varepsilon^l$ is the trivial $l$-plane bundle and $\gamma_{l+k}$ is the universal $(l+k)$-plane bundle; the jet bundle $J^\infty_0(\varepsilon^l,\gamma_{l+k})$ has as fibre the space $J^\infty_0(\R^l,\R^{l+k})$ of polynomial maps with $0$ constant term. Here we can consider subspaces $K(\R^l,\R^{l+k})\subset J^\infty_0(\R^l,\R^{l+k})$ invariant under the $\O(l+k)$-action induced by the one on $\R^{l+k}$ and such that the suspension of polynomial maps $J^\infty_0(\R^l,\R^{l+k})\to J^\infty_0(\R^{l+1},\R^{l+k+1})$ defined by $p\mapsto p\times\id_\R$ maps $K(\R^l,\R^{l+k})$ to $K(\R^{l+1},\R^{l+k+1})$; these subspaces yield a sequence of subbundles $K(\varepsilon^l,\gamma_{l+k})\subset J^\infty_0(\varepsilon^l,\gamma_{l+k})$. In the cases we are interested in, the subspace $\eta(\R^l,\R^{l+k})\subset J^\infty_0(\R^l,\R^{l+k})$ consisting of polynomials whose germ at $0$ has degeneracy type $\eta$ 
has the above properties which defines the subbundles $\eta(\varepsilon^l,\gamma_{l+k})\subset J^\infty_0(\varepsilon^l,\gamma_{l+k})$. Moreover, for each $l$ there is a cohomology class dual to the closure of $\eta(\varepsilon^l,\gamma_{l+k})$ denoted by $[\ol\eta(\varepsilon^l,\gamma_{l+k})]\in H^*(J^\infty_0(\varepsilon^l,\gamma_{l+k});\Z_2)\cong H^*(B\O(l+k);\Z_2)$ and the natural 
inclusion $B\O(l+k)\subset B\O(l+k+1)$ is such that the restriction $H^*(B\O(l+k+1);\Z_2)\to H^*(B\O(l+k);\Z_2)$ maps $[\ol\eta(\varepsilon^{l+1},\gamma_{l+k+1})]$ to $[\ol\eta(\varepsilon^l,\gamma_{l+k})]$.

\begin{defi}
  \label{defi:thpol}
  The \textit{Thom polynomial} of the singularity $\eta(k)$ is defined as the cohomology class $\Tp_{\eta(k)}\in H^*(B\O;\Z_2)\cong\Z_2[w_1,w_2,\ldots]$ that is the limit of the classes $[\ol\eta(\varepsilon^l,\gamma_{l+k})]\in H^*(B\O(l+k);\Z_2)$ as $l$ tends to $\infty$.
\end{defi}


For any map $f\colon M^n\to N^{n+k}$ the stable normal bundle $\nu_f$ is induced from $B\O$ by a map that we will (with a slight abuse of notation) also denote by $\nu_f\colon M\to B\O$. This map factors through a map $\nu_f(l)\colon M\to B\O(l+k)$ where $l$ is sufficiently large. Now the degeneracy locus $\eta(f)$ is the preimage of the subbundle $\nu_f(l)^*\eta(\varepsilon^l,\gamma_{l+k})\subset\kappa_f(l)^*J^\infty_0(\varepsilon^l,\gamma_{l+k})$ under the section defined by the jets of the germs of $f$, hence the Poincaré dual of $\ol\eta(f)$ is $\nu_f^*\Tp_{\eta(k)}$ which is the same as the polynomial obtained by substituting $w_i(\nu_f)$ in the variable $w_i$ of $\Tp_{\eta(k)}$ for all $i$.

\begin{defi}
  For a map $f\colon M\to N$, substituting the Stiefel--Whitney class $w_i(\nu_f)$ in the variable $w_i$ of $\Tp_{\eta(k)}$ for all $i$ will be called the \textit{evaluation} of the Thom polynomial on $f$ and denoted by $\Tp_{\eta(k)}(f)$.
\end{defi}

\begin{rmk}
  This definition allows us to evaluate Thom polynomials of $k$-codimensional singularities on non-$k$-codimensional maps too (this is used in corollaries \ref{crly:prim} and \ref{crly:twprim}) but if the map $f$ above is of codimension $k$, then the evaluation $\Tp_{\eta(k)}(f)$ is the cohomology class dual to $[\ol\eta(f)]\in H_*(M;\Z_2)$. If the map $f$ has codimension different from $k$, then one can consider an appropriate de-suspension (or suspension) $\tilde f$ of $f$ that has codimension $k$, and then the evaluation $\Tp_{\eta(k)}(f)$ gives the dual cohomology class of the closure of the $\eta$-locus $\ol\eta(\tilde f)\cap M$ if this intersection if transversal.
\end{rmk}

In the following we restrict ourselves to maps every singularity of which belongs to a fixed set of singularities $\tau$. Such a map will be called a \textit{$\tau$-map}. 
Then we can consider Thom polynomials modulo $\tau$, that is, Thom polynomials restricted to $\tau$-maps 
(see e.g. \cite{rim}) defined as follows.


\begin{defi}
  For any $l$ let $K_\tau(l)\subset J^\infty_0(\varepsilon^l,\gamma_{l+k})$ be the union of the subsets $\eta(\varepsilon^l,\gamma_{l+k})$ for all $\eta(k)\in\tau$ and let $K_\tau$ be the direct limit of the spaces $K_\tau(l)$ as $l$ tends to $\infty$. The space $K_\tau$ is called the \textit{Kazarian space} of the singularity set $\tau$ and the pullback of $\Tp_{\eta(k)}\in H^*(B\O;\Z_2)$ by the map $K_\tau\to B\O$ induced by the fibre maps $K_\tau(l)\to B\O(l+k)$ is called the \textit{Thom polynomial of $\eta(k)$ for $\tau$-maps} and denoted by $\Tp_{\eta(k)}|_\tau\in H^*(K_\tau;\Z_2)$.
\end{defi}

If $f\colon M^n\to N^{n+k}$ is a $\tau$-map, then the map $\nu_f\colon M\to B\O$ inducing the stable normal bundle $\nu_f$ factors through a map $\kappa_f\colon M\to K_\tau$ as shown in the diagram below, hence then we have $\Tp_{\eta(k)}(f)=\Tp_{\eta(k)}|_\tau(f)$ for all $\eta(k)\in\tau$.
$$\xymatrix{
  & K_\tau\ar[d]^\pi && \Tp_{\eta(k)}|_\tau\ar@{|->}[dl]_{\kappa_f^*} \\
  M\ar[ur]^{\kappa_f}\ar[r]_{\nu_f} & BO & \Tp_{\eta(k)}(f) & \Tp_{\eta(k)}\ar@{|->}[u]_{\pi^*}\ar@{|->}[l]^(.4){\nu_f^*}
}$$

\begin{ex}
  ~\vspace{-.5em}
  \begin{enumerate}
  \item If $\tau$ consists of all possible singularities of codimension-$k$ maps, we get $K_{\textrm{all}}\cong B\O$ and $\Tp_{\eta(k)}|_{\textrm{all}}=\Tp_{\eta(k)}$.
  \item For $\tau=\{\Sigma^0(k),\Sigma^1(k)\}$ we get $K_{\textrm{Morin}}$; computing the evaluation of $\Tp_{\Sigma^{1_r}(k)}$ on every Morin map (theorem \ref{thm:mor}) is the same as computing $\Tp_{\Sigma^{1_r}(k)}|_{\textrm{Morin}}$. 
  \item For $\tau=\{\Sigma^0(k),\Sigma^{1,0}(k),\Sigma^{1,1,0}(k)\}$ we get $K_{\textrm{cusp}}$.
  \end{enumerate}
\end{ex}


The Thom polynomials $\Tp_{\eta(k)}$ (and $\Tp_{\eta(k)}|_\tau$) considered so far give the classes represented by closures of $\eta$-loci in cohomology with $\Z_2$ coefficients. However, if the locus $\ol\eta(f)\subset M$ is cooriented for a map $f\colon M\to N$, then it also makes sense to talk about the cohomology class represented by it in $H^*(M;\Z)$. In many cases there are general conditions on singularities $\eta(k)$ which guarantee that the $\eta$-locus of any $k$-codimensional map is cooriented (see \cite{orsing}), but generally this is not sufficient for the existence of an analogue of the Thom polynomial $\Tp_{\eta(k)}$ (and $\Tp_{\eta(k)}|_\tau$) with $\Z$ coefficients; in fact a necessary and sufficient condition for the existence of an integral Thom polynomial for $\eta(k)$ is that it is a cocycle in \textit{Vassiliev's cochain complex} which is generated by cooriented singularities; see \cite{vas} (see also \cite{charsing}). If the integral Thom polynomial of $\eta(k)$ exists, then it is an element in the cohomology ring $H^*(B\O;\Z)$ which is isomorphic to $\Z[p_1,p_2,\ldots]\oplus\im\beta$ where $p_i$ denotes the $i$'th universal Pontryagin class and $\beta$ is the Bockstein map corresponding to the coefficient sequence $0\to\Z\to\Z\to\Z_2\to0$ (see \cite[problem 15-C]{charclass}). The restriction of the modulo $2$ reduction homomorphism is an isomorphism $\im\beta\to\im\Sq^1$ under which each element $w_{i_1}\ldots w_{i_r}\in\im\Sq^1$ corresponds to an element $v_{i_1,\ldots,i_r}\in\im\beta$. We note in particular that for $i$ odd we have $\Sq^1(w_iw_{i+1})=w_iw_{i+2}$, hence there is an element $v_{i,i+2}$ in $\im\beta$.

Since the existence of integral Thom polynomials is hard to decide in general (see \cite{inttp} for such computations) and the topic of the present paper is far from these techniques, we only state here the two cases we need and where this existence is easy to see. First note that since we consider here maps between oriented manifolds, the construction of Kazarian spaces can be repeated with $B\SO$ replacing $B\O$; the oriented analogue of $K_\tau$ will be denoted $K_\tau^\SO$, in other words, $K_\tau^\SO$ is the pullback of $K_\tau$ by the double cover $B\SO\to B\O$. (Note that the properties of the cohomology ring of $B\O$ recalled above also hold for that of $B\SO$.)

\begin{defi}
  If the singularity $\eta(k)$ is a cocycle in the restricted Vassiliev complex that is generated only by the cooriented singularities of oriented maps belonging to $\tau$, then there is a corresponding cohomology class in the ring $H^*(K_\tau^\SO;\Z)$ called the \textit{integral Thom polynomial of $\eta(k)$ for $\tau$-maps} and denoted by $\Tp^\SO_{\eta(k)}|_\tau$.
\end{defi}

Integral Thom polynomials are analogous to Thom polynomials in $\Z_2$-cohomology. Namely for any $\tau$-map $f\colon M^n\to N^{n+k}$ between oriented manifolds, if $\eta(k)\in\tau$ satisfies the condition of the above definition, then the Poincaré dual of the closure $\ol\eta(f)$ in $H^*(M;\Z)$ is the class $\kappa_f^*\Tp_{\eta(k)}^\SO|_\tau$ which is obtained by substituting $p_i(\nu_f)$ and $v_{i_1,\ldots,i_r}(\nu_f)$ in the variables $p_i$ and $v_{i_1,\ldots,i_r}$ of this Thom polynomial respectively.

The first case in which we will consider integral Thom polynomials is that of Morin maps $f\colon M^n\to N^{n+k}$ where $M$ and $N$ are oriented and $k$ is odd; in this case the only possible singular loci are the $\Sigma^{1_r}(f)$ and among them the cooriented ones are precisely those where $r$ is even (see \cite{orsing}). This means that in the Vassiliev complex for only odd codimensional Morin maps between oriented manifolds, all differentials vanish, hence $\Sigma^{1_r}(k)$ admits a Thom polynomial $\Tp^\SO_{\Sigma^{1_r}(k)}|_{\textrm{Morin}}$ in this case.



Secondly we consider generic maps $f\colon M^n\to N^{n+l}$ where $M$ and $N$ are oriented; in this case the Thom--Boardman stratum $\Sigma^r(f)$ is cooriented precisely if $l$ is even (see \cite{orsing}) and then \cite{calcsk} shows that the Thom polynomial $\Tp^\SO_{\Sigma^r(l)}$ of $\Sigma^r(l)$ exists precisely if $r$ is also even.
(Note that in the cases we are interested in we have $l=k-r+1$ and so if $r$ is even, then $k$ is odd precisely if $l$ is even.)





\section{Proofs}
\label{sec:prf}

\begin{prf}[of theorem \ref{thm:cusp}]
  For any map $f\colon M^n\to N^{n+k}$ the composition $f\circ\pr_M\colon M\times\R\to N$ (which is not a generic map) is such that we have $\Sigma^2(f\circ\pr_M)\cap M\times\{0\}=\Sigma^1(f)$. Our aim is to show that in the case when $f$ is a cusp map, any generic perturbation of this composition for which the submanifold $M\times\{0\}$ is transverse to the singular strata is such that its $\Sigma^2$-locus intersected with $M\times\{0\}$ is cobordant to the $\Sigma^{1_2}$-locus of $f$. Below we construct a map $\tilde f\colon M\times\R\to N$ such that we will have $\Sigma^2(\tilde f)\cap M\times\{0\}=\Sigma^{1_2}(f)$ (not just up to cobordism) where $\tilde f$ is a perturbation of $f\circ\pr_M$ but it is (in general) not generic. However, the $1$-jet $j^1\tilde f\colon M\times\R\to J^1(M\times\R,N)$ will be transverse to the subspace $\Sigma^2(M\times\R,N)$ of $1$-jets of corank $2$. This is sufficient to prove our statement since both being a $\Sigma^2$-map (i.e. having no singularities of type $\Sigma^3$) and having a $1$-jet extension transverse to $\Sigma^2(M\times\R,N)$ are open conditions, 
  hence any generic perturbation $\hat f$ of $f\circ\pr_M$ is homotopic to $\tilde f$ through $\Sigma^2$-maps whose $1$-jets are transverse to $\Sigma^2(M\times\R,N)$. Now if we also suppose transverse intersection of the singular strata with $M\times\{0\}$ (and with $M\times\{0\}\times[0,1]$ for the homotopy), we obtain that $\Sigma^2(\hat f)\cap M\times\{0\}$ is cobordant to $\Sigma^2(\tilde f)\cap M\times\{0\}$ as submanifolds of $M$ which is what we claimed.


  For any section $\sigma$ of the vector bundle $f^*TN$ over $M$ we can perturb the map $f\circ\pr_M$ by $\sigma$ by mapping the point $(p,t)\in M\times\R$ to the point $\exp_N(t\sigma(p))$ where $\exp_N$ is the exponential map of $N$ corresponding to a Riemannian metric; 
  this may only be well-defined on a small neighbourhood of $M\times\{0\}$ but that can be identified with $M\times\R$. The map given by this formula will be denoted $\tilde f_\sigma\colon M\times\R\to N$. 
  Observe that here $\Sigma^2(\tilde f_\sigma)\cap M\times\{0\}$ is the zero locus $\sigma^{-1}(0)$ intersected with $\Sigma(f)$, hence our goal is now to find a section $\sigma\colon M\to f^*TN$ for which $\sigma^{-1}(0)\cap\Sigma(f)$ is precisely $\Sigma^{1_2}(f)$. We also note that generally the map $\tilde f_\sigma$ is not finite-to-one, thus not generic.

  Porteous \cite{gtp} showed that for any map $f$ there is a second intrinsic derivative
  $$d^2f\colon\ker df\otimes\ker df\to\coker df$$
  of $f$ given by the 
  usual second derivative 
  in the direction of the subspace $\ker df<TM$, and its value is well-defined modulo the subspace $\im df<f^*TN$, i.e. it is independent of the choice of local coordinates. Now $\ker df|_{\Sigma(f)}$ is a line bundle and so its tensor square is trivial, hence the second derivative $d^2f|_{\Sigma(f)}$ defines a section of the bundle $\coker df|_{\Sigma(f)}$ which we identify with the orthogonal complement of $\im df|_{\Sigma(f)}$ using a Riemannian metric. This gives us a section of the restriction of $f^*TN$ to $\Sigma(f)$. Since the second derivative vanishes precisely at $\Sigma^{1_2}(f)$, so does this section, hence if we denote by $\sigma$ an arbitrary extension of it to the whole $M$, we can define $\tilde f$ as $\tilde f_\sigma$.

  It is left for us to show that the $1$-jet $j^1\tilde f$ is transverse to $\Sigma^2(M\times\R,N)$. This is done by writing the map $f$ in local coordinates around the points of $\Sigma^{1_2}(f)\times\{0\}$ (i.e. writing up the normal form of cusp singularity; see \cite{mor}), computing from it the map $\tilde f$ in these coordinates, then proving that the derivative $\R^{n+1}\to\Hom(\R^{n+1},\R^{n+k})$ of this local map is transverse to $\Sigma^2(\R^{n+1},\R^{n+k})$. A way of this is showing that at each $\Sigma^2$-point its derivative projected to the normal space of $\Sigma^2(\R^{n+1},\R^{n+k})$ is surjective using that the normal space of $\Sigma^2(\R^{n+1},\R^{n+k})$ at a point $\varphi\colon\R^{n+1}\to\R^{n+k}$ is $\Hom(\ker\varphi,\coker\varphi)$ (see \cite{calcord2}). This computation is lengthy but straightforward and so we omit the details here but the interested reader can find them in the \hyperref[sec:app]{appendix}. This concludes the proof.
\end{prf}

\begin{prf}[of theorem \ref{thm:prim}]
  We will prove the following:

  \begin{thm'}{thm:prim}
    If $f\colon M^n\to N^{n+k}$ is a prim map where $M$ is closed, $k\ge r-1$ and $N$ is of the form $Q^{n+k-r+1}\times P^{r-1}$ such that the composition $g:=\pr_Q\circ f$ is generic and the pullback bundle $f^*\pr_P^*TP$ is trivial, then the submanifolds $\ol\Sigma^{1_r}(f)\subset M$ and $\Sigma^r(g)\subset M$ are cobordant as embedded submanifolds of $M$. This cobordism is oriented if $M$ and $N$ are oriented and the codimension $k$ is odd.
  \end{thm'}

  This implies theorem \ref{thm:prim} since if $f\colon M\to N$ is any prim map and $\tilde f\colon M\times\R^{r-1}\to N$ is its de-suspension as in theorem \ref{thm:prim} (we may also assume $\tilde f|_{M\times\{0\}}=f$), then we can use theorem \hyperref[thm:prim']{\ref*{thm:prim}'} for maps $\hat f\colon M\times T^{r-1}\to N\times T^{r-1}$ and $\hat g\colon M\times T^{r-1}\to N$ where $T^{r-1}$ is the $(r-1)$-torus, $\hat f$ is a perturbation of $f\times\id_{T^{r-1}}$ and $\hat g$ is $\pr_N\circ\hat f$ such that on a small neighbourhood of $M\times\{0\}$ the maps $\hat g$ and $\tilde f$ coincide. We can indeed apply theorem \hyperref[thm:prim']{\ref*{thm:prim}'} for $\hat f$ and $\hat g$ since $\hat f$ is a prim map and $T^{r-1}$ is parallelisable implying that the pullback of its tangent bundle is trivial. Thus we get that $\ol\Sigma^{1_r}(\hat f)$ and $\Sigma^r(\hat g)$ are cobordant submanifolds of $M\times T^{r-1}$ and by intersecting them with $M\times\{0\}$ (making the cobordism transverse to it) we get that $\ol\Sigma^{1_r}(f)$ and $\Sigma^r(\tilde f)\cap M\times\{0\}$ are cobordant in $M$.
  
  Now a trivial observation is that if $\xi$ is a smooth vector bundle over a closed manifold $M$ and $\sigma_1,\ldots,\sigma_r$ and $\tau_1,\ldots,\tau_r$ are two collections of sections, both in general position, 
  then the zero sets $\bigcap\limits_{i=1}^r\sigma_i^{-1}(0)$ and $\bigcap\limits_{i=1}^r\tau_i^{-1}(0)$ are cobordant submanifolds of $M$. This cobordism is oriented if $M$ and $\xi$ are oriented and $\xi$ is of even rank. Theorem \hyperref[thm:prim']{\ref*{thm:prim}'} will follow if we substitute for $\xi$ the normal bundle $\nu_j$ of 
  an immersion lift $j\colon M\imto N\times\R$ of the prim map $f$ and use the following two lemmas:

  \begin{lemma}
    \label{lemma:1r}
    If $f\colon M\to N$ is a prim map and $j\colon M\imto N\times\R$ is its immersion lift, then for any $r$ there are sections $\sigma_1,\ldots,\sigma_r$ (in general position) of $\nu_j$ (the normal bundle of $j$) such that $\ol\Sigma^{1_r}(f)$ is the set $\bigcap\limits_{i=1}^r\sigma_i^{-1}(0)$.
  \end{lemma}
  
  \begin{lemma}
    \label{lemma:r}
    If $j\colon M\imto Q\times P^r$ is an immersion where the composition $g:=\pr_Q\circ j$ is generic and the pullback bundle $j^*\pr_P^*TP$ is trivial, then there are sections $\tau_1,\ldots,\tau_r$ (in general position) of $\nu_j$ such that $\Sigma^r(g)$ is the set $\bigcap\limits_{i=1}^r\tau_i^{-1}(0)$.\footnote{Note that we apply lemma \ref{lemma:r} in the proof of theorem \hyperref[thm:prim']{\ref*{thm:prim}'} with the (strange) substitution $P^r=P^{r-1}\times\R$.}
  \end{lemma}

  To see that lemma \ref{lemma:r} holds, note that
  $$j^*T(Q\times P)\cong j^*\textstyle\pr_Q^*TQ\oplus\varepsilon^r\cong TM\oplus\nu_j$$
  (where $\varepsilon^r$ again denotes the trivial $r$-plane bundle) and the standard coordinate vectors in $\varepsilon^r$ yield sections $\varepsilon_1,\ldots,\varepsilon_r$ of $TM\oplus\nu_j$. Then 
  letting $\tau_i$ be the projection of $\varepsilon_i$ to the normal bundle $\nu_j$ gives the sections of $\nu_j$ with the desired property since $\Sigma^r(g)$ is the set of points $p\in M$ where the whole summand $\varepsilon^r|_p<j^*(T(Q\times P))|_p$ is in the tangent space $T_pM$. (Throughout the paper we tacitly identify the normal bundle of a submanifold with the orthogonal complement of its tangent bundle and suppose that the projection to the normal bundle is the orthogonal projection.)

  Lemma \ref{lemma:1r} has been proved in \cite[proof of lemma 2]{singprim}, for the reader's convenience we briefly recall the idea here. We apply lemma \ref{lemma:r} with $Q=N$ and $P=\R$ to obtain $\sigma_1:=\tau_1=\pr_{\nu_j}\circ\varepsilon$ where $\varepsilon:=\varepsilon_1$ is the constant vector field on $N\times\R$ defined by the positive basis of $\R$, hence we have $\sigma_1^{-1}(0)=\Sigma(f)$. Now $\Sigma(f)$ is a submanifold of $M$ of codimension $k+1$; we denote its normal bundle in $M$ by $\nu_1$ and observe that the tangent spaces of $\sigma_1$ at the points of $\Sigma(f)$ form the graph of a pointwise linear isomorphism $\iota_1\colon\nu_1\to\nu_j|_{\Sigma(f)}$. The composition $\iota_1\circ\pr_{\nu_1}\circ\varepsilon$ is a section of $\nu_j|_{\Sigma(f)}$ which vanishes precisely at $\Sigma^{1,1}(f)$, hence if we define $\sigma_2$ as an arbitrary extension of this section to the whole $M$, then we obtain $\sigma_1^{-1}(0)\cap\sigma_2^{-1}(0)=\Sigma^{1,1}(f)$. By further applying the same method, we obtain the sections $\sigma_3,\ldots,\sigma_r$ with the desired property.
\end{prf}

\begin{prf}[of theorem \ref{thm:twprim}]
  We will prove the following, which implies theorem \ref{thm:twprim} analogously to how theorem \hyperref[thm:prim']{\ref*{thm:prim}'} implied theorem \ref{thm:prim}.

  \begin{thm'}{thm:twprim}
    If $f\colon M^n\to N^{n+k}$ is a twisted prim map where $M$ is closed and $N$ is of the form $Q^{n+k-r+1}\times P^{r-1}$ such that the composition $g:=\pr_Q\circ f$ is generic and the pullback bundle $f^*\pr_P^*TP$ is trivial, further, $M$ and $N$ are oriented, $k$ is odd and $r$ is even, then the homology classes $2[\ol\Sigma^{1_r}(f)],2[\Sigma^r(g)]\in H_*(M;\Z)$ coincide.
  \end{thm'}
  
  To prove theorem \hyperref[thm:twprim']{\ref*{thm:twprim}'} let $\ell\to M$ be a line bundle with $\ell|_{\Sigma(f)}\cong\ker df$, let $q\colon\tilde M\to M$ be its $S^0$-bundle and put $\tilde f:=f\circ q$ and $\tilde g:=g\circ q$. Now we have $\ker d\tilde f=q^*\ker df\cong q^*\ell|_{\Sigma(f)}$ which is trivial, hence the map $\tilde f\colon\tilde M\to N$ is prim and so by theorem \hyperref[thm:prim']{\ref*{thm:prim}'} the submanifolds $\ol\Sigma^{1_r}(\tilde f)$ and $\Sigma^r(\tilde g)$ are cobordant in $\tilde M$. In particular we have
  $$[\ol\Sigma^{1_r}(\tilde f)]=[\Sigma^r(\tilde g)]\in H_*(\tilde M;\Z).$$
  
  The map $q$ is a double cover, hence for each component $C$ of $\ol\Sigma^{1_r}(f)$ or $\Sigma^r(g)$, the homology class $q_*[q^{-1}C]$ is a multiple of $[C]$ by either $0$ or $2$. Since we also have $\ol\Sigma^{1_r}(\tilde f)=q^{-1}\ol\Sigma^{1_r}(f)$ and $\Sigma^r(\tilde g)=q^{-1}\Sigma^r(g)$, it is sufficient to prove that $q_*[q^{-1}C]$ is always $2[C]$, or, in other words, the orientation gained on the simplices of $\ol\Sigma^{1_r}(\tilde f)$ (resp. $\Sigma^r(\tilde g)$) from the map $\tilde f$ (resp. $\tilde g$) is the same as the one gained from lifting the orientation of the simplices of $\ol\Sigma^{1_r}(f)$ (resp. $\Sigma^r(g)$).

  Now the orientation gained from $\tilde f$ (resp. $\tilde g$) is the same as the natural orientation of $\bigcap\limits_{i=1}^r\sigma_i^{-1}(0)$ for $r$ generic sections $\sigma_1,\ldots,\sigma_r$ of $\nu_{j}$ where $j\colon\tilde M\to N\times\R$ is an immersion lift of $\tilde f$. These sections are (resp. the first of these sections is) obtained from the standard unit vector field on $N\times\R$ whose pullback to $\Sigma(\tilde f)$ is $\ker d\tilde f$ which extends to $\tilde M$ as $q^*\ell$. This means that for each point $p\in M$, the two fibres of $q^*\ell|_{q^{-1}p}$ are trivialised by vectors whose projections to $\ell_p$ are opposite. Thus for a small neighbourhood $U\subset M$ of the point $p\in M$ the two local sections are each other's opposite, i.e. if the two components of $q^{-1}U$ are $\tilde U_1,\tilde U_2\subset\tilde M$ and $\varphi\colon\nu_{j}|_{\tilde U_1}\to\nu_{j}|_{\tilde U_2}$ is the natural isomorphism, the section $\varphi\circ\sigma_i|_{\tilde U_1}$ is $-\sigma_i|_{\tilde U_2}$ for all $i$ (resp. for $i=1$). But since the codimension $k$ is odd, the rank of $\nu_{j}$ is even, hence the multiplication by $-1$ does not change the (local) orientation of the zero locus of the section. This is what we wanted.
\end{prf}


The proofs of the corollaries we will be based on the following: 

\begin{lemma}
  \label{lemma:eq}
  Let $f\colon M^n\to N^{n+k}$ be a generic $\tau$-map where $\tau$ is a set of singularities, 
  let $\tilde f\colon M\times\R^l\to N$ be a de-suspension of $f$ and suppose that $\eta(k)$ and $\vartheta(k-l)$ are singularities for which the Thom polynomials $\Tp^\circ_{\eta(k)}|_\tau$ and $\Tp^\circ_{\vartheta(k-l)}$ exist where $\circ$ stands for $\SO$ if $M$ and $N$ are oriented and nothing otherwise. If the singular loci $\eta(f)$ and $\vartheta(\tilde f)$ are such that for some numbers $a,b$ we have $a[\ol\eta(f)]=b[\ol\vartheta(\tilde f)\cap M\times\{0\}]$ in $H_*(M;\Z)$ in the oriented case and in $H_*(M;\Z_2)$ otherwise, then we have
  $$a\Tp^\circ_{\eta(k)}(f)=b\Tp^\circ_{\vartheta(k-l)}(f).
  $$
\end{lemma}

\begin{prf}
  This is trivial since the homology classes $[\ol\eta(f)]$ and $[\ol\vartheta(\tilde f)\cap M\times\{0\}]$ are dual to the cohomology classes $\Tp^\circ_{\eta(k)}(f)$ and $\Tp^\circ_{\vartheta(k-l)}(f)$ respectively.
\end{prf}



\begin{prf}[of corollary \ref{crly:cusp}]
  Theorem \ref{thm:cusp} and lemma \ref{lemma:eq} imply that for all \textit{cusp} maps $f\colon M^n\to N^{n+k}$ we have $\Tp_{\Sigma^{1_2}(k)}(f)=\Tp_{\Sigma^2(k-1)}(f)$ and if $k$ is odd and $M$ and $N$ are oriented, then the same with integral Thom polynomials. However, this already proves that the same holds for \textit{any generic} map $f\colon M^n\to N^{n+k}$. Indeed, for generic maps new singularities can occur compared to the case of cusp maps but these new singularities have large codimensions (greater than $2k+3$) hence they do not interfere with cohomologies in the codimension of the cusp singularity $\Sigma^{1_2}(k)$ (which is $2k+2$). This is because for any $r$ the codimension of the singularity $\Sigma^{1_r}(k)$ 
  is $r(k+1)$ and the codimension of $\Sigma^r(k)$ is $r(k+r)$. 
  An alternative construction of Kazarian spaces gives them by gluing together disk bundles of vector bundles corresponding to singularities where the rank of each vector bundle is the codimension of the singularity it corresponds to; see e.g. \cite{calctp}, \cite{hosszu}. This implies that the embedding of the Kazarian space $K^\circ_{\textrm{cusp}}$ of the singularity set $\{\Sigma^0(k),\Sigma^{1,0}(k),\Sigma^{1,1,0}(k)\}$ (where $\circ$ again denotes $\SO$ or nothing in the oriented and non-oriented cases respectively) into $K^\circ_{\textrm{all}}=:B$ (which is $B\SO$ and $B\O$ in the oriented and non-oriented cases respectively) is such that the additional cells of $B$ 
  are all of dimension at least $2k+4$. Thus we have $H^i(B)\cong H^i(K^\circ_{\textrm{cusp}})$ for $i\le2k+2$ with any coefficient group which means that $\Tp^\circ_{\Sigma^{1_2}(k)}|_{\textrm{cusp}}\in H^{2k+2}(K^\circ_{\textrm{cusp}})$ corresponds isomorphically to the element $\Tp^\circ_{\Sigma^{1_2}(k)}\in H^{2k+2}(B)$. (Note that this also means that the integral Thom polynomial $\Tp^\SO_{\Sigma^{1_2}(k)}$ exists even for non-Morin maps.)

  But the fact that for all maps $f\colon M^n\to N^{n+k}$ we have $\Tp^\circ_{\Sigma^{1_2}(k)}(f)=\Tp^\circ_{\Sigma^2(k-1)}(f)$ 
  is equivalent to saying that for any manifold $M$ and any map $\nu\colon M\to B$ the pullbacks $\nu^*\Tp^\circ_{\Sigma^{1_2}(k)}$ and $\nu^*\Tp^\circ_{\Sigma^2(k-1)}$ coincide. Since the space $B$ is the direct limit of the sequence of Grassmannian manifolds, this can only be if the cohomology classes $\Tp^\circ_{\Sigma^{1_2}(k)}$ and $\Tp^\circ_{\Sigma^2(k-1)}$ are also the same.
\end{prf}

  

\begin{prf}[of corollary \ref{crly:prim}]
  From theorem \ref{thm:prim} and lemma \ref{lemma:eq} we obtain for any prim map $f\colon M^n\to N^{n+k}$ the equality $\Tp_{\Sigma^{1_r}(k)}(f)=\Tp_{\Sigma^r(k-r+1)}(f)$ and if $k$ is odd, $r$ is even and $M$ and $N$ are oriented, then the same with integer Thom polynomials. The fact that this cohomology class equals $w_{k+1}(\nu_j)^r$ in the non-oriented and $p_{\frac{k+1}2}(\nu_j)^{\frac r2}$ in the oriented case follows from lemma \ref{lemma:1r} (note that $\nu_j\cong\nu_f\oplus\varepsilon^1$). Indeed, the Poincaré dual of the zero set $\sigma_i^{-1}(0)$ of a generic section $\sigma_i$ of $\nu_j$ is the Euler class $e(\nu_j)$, hence that of the intersection $\bigcap\limits_{i=1}^r\sigma_i^{-1}(0)$ is $e(\nu_j)^r$ which is $w_{k+1}(\nu_j)^r$ modulo $2$ and $p_{\frac{k+1}2}(\nu_j)^{\frac r2}$ if $\nu_j$ is oriented and $k$ is odd and $r$ is even.
\end{prf}

\begin{prf}[of corollary \ref{crly:twprim}]
  This is immediate from theorem \ref{thm:twprim} and lemma \ref{lemma:eq}.
\end{prf}

It is left for us to prove theorem \ref{thm:mor} for which we will need:

\begin{lemma}[Kazarian \cite{morchar}]
  \label{lemma:nu1}
  If $f\colon M^n\to N^{n+k}$ is a Morin map, then the normal bundle of the singular set $\Sigma(f)$ in $M$ is stably isomorphic to $\ker df\otimes\nu_f|_{\Sigma(f)}$ and if $i\colon\Sigma(f)\into M$ is the embedding, then for any $r\ge0$ the pushforward $i_!\left(w_1(\ker df)^r\right)$ is $w_{k+r+1}(\nu_f)$.
\end{lemma}

\begin{prf}
  Let $\nu$ be the normal bundle of $\Sigma(f)\subset M$ and let $q\colon PTM\to M$ be the projectivisation of the tangent bundle of $M$ and $\gamma_{PTM}$ the tautological line bundle over $PTM$. Consider the restriction of the bundle $q$ over the singular set, i.e. $q^{-1}\Sigma(f)\to\Sigma(f)$. Associating to each point $p\in\Sigma(f)$ the line $\ker df_p<TM$ defines a section $s\colon\Sigma(f)\to q^{-1}\Sigma(f)$ of this bundle. Note that $\ker df\cong s^*\gamma_{PTM}$. If $\nu_s$ denotes the normal bundle of $s(\Sigma(f))$ in $PTM$, then we have $\nu_s\cong q^*\nu\oplus\ker dq|_{s(\Sigma(f))}$, hence the vector bundle $q^*\nu$ is stably isomorphic to $\nu_s\ominus\ker dq|_{s(\Sigma(f))}$. On the other hand the differential $df$ defines a map of the line bundle $\gamma_{PTM}$ to $q^*f^*TN$, i.e. a section of the bundle $\Hom(\gamma_{PTM},q^*f^*TN)\cong\gamma_{PTM}\otimes q^*f^*TN$. The zero locus of this section is precisely $s(\Sigma(f))$, hence the normal bundle $\nu_s$ is isomorphic to the restriction of $\gamma_{PTM}\otimes q^*f^*TN$ to $s(\Sigma(f))$. Finally, observe that the tangent space of a projective space is isomorphic to the space of homomorphisms of the tautological line bundle to its orthogonal complement. Applying this observation fibrewise to the bundle $\ker dq\to PTM$ we obtain that $\ker dq$ is $\gamma_{PTM}\otimes(q^*TM/\gamma_{PTM})$ which is stably isomorphic to $\gamma_{PTM}\otimes q^*TM\ominus\varepsilon^1$. Putting all of the above observations together we obtain
  \begin{alignat*}2
    \nu&\cong s^*q^*\nu\cong s^*(\nu_s\ominus\ker dq|_{s(\Sigma(f))})\cong s^*(\gamma_{PTM}\otimes q^*f^*TN\ominus\ker dq)\cong \\
       &\cong s^*(\gamma_{PTM}\otimes q^*f^*TN\ominus(\gamma_{PTM}\otimes q^*TM\ominus\varepsilon^1))\cong s^*(\gamma_{PTM}\otimes q^*\nu_f\oplus\varepsilon^1)\cong \\
       &\cong\ker df\otimes\nu_f|_{\Sigma(f)}\oplus\varepsilon^1.
  \end{alignat*}


  To compute the image of $w_1(\ker df)^r$ under $i_!$, note that the pushforward $s_!\left(w_1(\ker df)^r\right)$ is $w_1(\gamma_{PTM})^r$ times the cohomology class represented by $s(\Sigma(f))$ in $PTM$. This latter cohomology class is $w_{n+k}(\gamma_{PTM}\otimes q^*f^*TN)$ since the submanifold representing it is the zero locus of a section of this bundle. We also know the pushforward image $q_!\left(w_1(\gamma_{PTM})^r\right)$, it is the normal Stiefel--Whitney class $\ol w_{r-n+1}(M)$ (see e.g. \cite{damonflag}), hence we have
  \begin{alignat*}2
    i_!\left(w_1(\ker df)^r\right)&=q_!s_!\left(w_1(\ker df)^r\right)=q_!\left(w_1(\gamma_{PTM})^r\sum_{i=0}^{n+k}w_1(\gamma_{PTM})^iw_{n+k-i}(q^*f^*TN)\right)=\\
                                  &=q_!\left(\sum_{i=r}^{n+k+r}w_1(\gamma_{PTM})^iw_{n+k+r-i}(q^*f^*TN)\right)=\\
                                  &=\sum_{i=r-n+1}^{k+r+1}\ol w_i(TM)w_{k+r+1-i}(f^*TN)=w_{k+r+1}(\nu_f).
  \end{alignat*}
  This is what we wanted.
\end{prf}

\begin{prf}[of theorem \ref{thm:mor}]
  For simplicity of notation, let $\ell$ denote the line bundle $\ker df$ (which is now only defined over $\Sigma(f)$) and denote for all $r\ge1$ by $\nu_r$ the normal bundle of the singular locus $\ol\Sigma^{1_r}(f)$ in $\ol\Sigma^{1_{r-1}}(f)$ (noting that $\ol\Sigma^{1_0}(f)=M$) and by $i_r$ the embedding $\ol\Sigma^{1_r}(f)\into M$. The following will be analogous to the proof of lemma \ref{lemma:1r}.

  Taking for all points $p\in\Sigma(f)$ the projection of the line $\ell|_p<T_pM$ to the normal space $\nu_1|_p$ defines a section $\sigma_2$ of the bundle $\Hom(\ell,\nu_1)\cong\ell\otimes\nu_1$. Observe that the vanishing locus $\sigma_2^{-1}(0)$ is $\Sigma^{1,1}(f)$ which is a submanifold in $\Sigma(f)$ of codimension $k+1$ with normal bundle $\nu_2$ and note that the tangent spaces of $\sigma_2$ at the points of $\Sigma^{1,1}(f)$ form the graph of a pointwise isomorphism $\iota_2\colon\nu_2\to\ell\otimes\nu_1|_{\Sigma^{1,1}(f)}$. The projection of $\ell|_{\Sigma^{1,1}(f)}$ to $\nu_2$ composed with $\iota_2$ gives a section of the bundle $\Hom(\ell,\ell\otimes\nu_1)|_{\Sigma^{1,1}(f)}\cong\nu_1|_{\Sigma^{1,1}(f)}$ which vanishes precisely at $\Sigma^{1,1,1}(f)$, hence if $\sigma_3$ is an arbitrary extension of this section to the whole $\Sigma(f)$, then we get that $\Sigma^{1,1,1}(f)$ is the intersection $\sigma_2^{-1}(0)\cap\sigma_3^{-1}(0)$ where $\sigma_2$ is a section of $\ell\otimes\nu_1$ and $\sigma_3$ is a section of $\nu_1$. Further iterating this method produces the sequence $\sigma_2,\sigma_3,\ldots$ where each term of even (resp. odd) index is a section of the bundle $\ell\otimes\nu_1$ (resp. $\nu_1$) such that for any $r$ we have $\ol\Sigma^{1_r}(f)=\bigcap\limits_{i=1}^r\sigma_i^{-1}(0)$.

  Thus the dual cohomology class represented by $\ol\Sigma^{1_r}(f)$ in the submanifold $\Sigma(f)\subset M$ is $e(\ell\otimes\nu_1)^{\lfloor \frac r2 \rfloor}e(\nu_1)^{\lceil \frac r2 \rceil-1}$ where by lemma \ref{lemma:nu1} we have $\nu_1\cong\ell\otimes\nu_f\oplus\varepsilon^1$, hence $\ell\otimes\nu_1\cong\nu_f\oplus\ell$ (not denoting the restriction of $\nu_f$ to $\Sigma(f)$ for simplicity). 
  Note that over $\Sigma(f)$ the expression $\nu_f\oplus\ell$ can be considered as a genuine rank-$(k+1)$ vector bundle, not just a stable one. Hence the modulo $2$ Euler classes $e(\ell\otimes\nu_1)$ and $e(\nu_1)$ are $w_{k+1}(\nu_f\oplus\ell)=w_{k+1}(\nu_f)+w_1(\ell)w_k(\nu_f)$ and $w_{k+1}(\ell\otimes\nu_f)$ respectively. We claim that the latter coincides with $w_{k+1}(\nu_f)$. Indeed, another expansion of $w_{k+1}(\nu_f\oplus\ell)$ is
  \begin{alignat*}2
    w_{k+1}\left(\ell\otimes(\ell\otimes\nu_f\oplus\varepsilon^1)\right)&=w_{k+1}(\ell\otimes\nu_f)+w_1(\ell)\sum_{i=0}^kw_1(\ell)^iw_{k-i}(\ell\otimes\nu_f)=\\
                                                                        &=w_{k+1}(\ell\otimes\nu_f)+w_1(\ell)w_k(\nu_f),
  \end{alignat*}
  hence $w_{k+1}(\ell\otimes\nu_f)$ is $w_{k+1}(\nu_f\oplus\ell)+w_1(\ell)w_k(\nu_f)=w_{k+1}(\nu_f)$. Thus we have
  $$\Tp_{\Sigma^{1_r}(k)}(f)=(i_1)_!\left(\left(w_{k+1}(\nu_f)+w_1(\ell)w_k(\nu_f)\right)^{\lfloor \frac r2 \rfloor}w_{k+1}(\nu_f)^{\lceil \frac r2 \rceil-1}\right)$$
  which is the pushforward by $i_1$ of $\left(w_{k+1}(\nu_f)^2+w_1(\ell)w_{k+1}(\nu_f)w_k(\nu_f)\right)^{\frac{r-1}2}$ if $r$ is odd and of $(w_{k+1}(\nu_f)+w_1(\ell)w_k(\nu_f))\left(w_{k+1}(\nu_f)^2+w_1(\ell)w_{k+1}(\nu_f)w_k(\nu_f)\right)^{\frac r2-1}$ if $r$ is even. Note that since the stable vector bundle $\nu_f\oplus\ell$ has as a representative the rank-$(k+1)$ bundle $\nu_1$, we have $0=w_{k+2}(\nu_f\oplus\ell)=w_{k+2}(\nu_f)\oplus w_1(\ell)w_{k+1}(\nu_f)$, i.e. $w_1(\ell)w_{k+1}(\nu_f)=w_{k+2}(\nu_f)$. Hence by lemma \ref{lemma:nu1} we get that the cohomology class $\Tp_{\Sigma^{1_r}(k)}(f)$ is $w_{k+1}(\nu_f)\left(w_{k+1}(\nu_f)^2+w_{k+2}(\nu_f)w_k(\nu_f)\right)^{\frac{r-1}2}$ if $r$ is odd and $\left(w_{k+1}(\nu_f)^2+w_{k+2}(\nu_f)w_k(\nu_f)\right)^{\frac r2}$ if $r$ is even.

  Let us now turn to the case where $k$ is odd, $r$ is even and $M$ and $N$ are oriented. Then $\Sigma(f)\subset M$ is generally not a cooriented submanifold but $\Sigma^{1,1}(f)\subset M$ is. The dual cohomology class represented by $\ol\Sigma^{1_r}(f)$ in $\Sigma^{1,1}(f)$ is well-defined with $\Z$ coefficients; it is the Euler class power $e(\ell\otimes\nu_1\oplus\nu_1)^{\frac r2-1}$ (again we suppress the restriction of $\nu_1$ and $\ell$ to $\Sigma^{1,1}(f)$ for notational simplicity). Note that it is sufficient to prove that the double $2\Tp^\SO_{\Sigma^{1_r}(k)}(f)$ of the cohomology class represented by $\ol\Sigma^{1_r}(f)$ is $2p_{\frac{k+1}2}(\nu_f)^{\frac r2}$. Indeed, if this holds, then the difference of $\Tp^\SO_{\Sigma^{1_r}(k)}(f)$ and $p_{\frac{k+1}2}(\nu_f)^{\frac r2}$ is an element of second order which is determined by the modulo $2$ restrictions of this Thom polynomial and Pontryagin class; these are $\Tp_{\Sigma^{1_r}(k)}(f)$ and $w_{k+1}(\nu_f)^r$ respectively which implies our formula. Thus we first want to determine the double of an Euler class for which we can use:

  \begin{lemma}
    \label{lemma:eu}
    If $q\colon\tilde M\to M$ is a double cover of an oriented manifold $M$ and $\xi$ is an oriented vector bundle over $M$, then we have $q_!q^*e(\xi)=2e(\xi)$.
  \end{lemma}

  Lemma \ref{lemma:eu} simply follows from the fact that if $\sigma$ is a generic section of $\xi$, then $q^*e(\xi)$ is the Poincaré dual of the zero locus of the section $\tilde\sigma$ of $q^*\xi$ induced from $\sigma$ and for each simplex $\Delta$ of $\sigma^{-1}(0)$, the preimage $q^{-1}\Delta$ is two simplices of $\tilde\sigma^{-1}(0)$ with the same orientation, hence we have $q_*q^{-1}\Delta=2\Delta$.

  We use lemma \ref{lemma:eu} with $\xi=\left(\frac r2-1\right)(\ell\otimes\nu_1\oplus\nu_1)$ and $q\colon\tilde M\to M$ the $S^0$-bundle of $\ell$. This bundle $\xi$ is indeed orientable since the rank of $\nu_1$ is $k+1$ which is even, hence we have $w_1(\ell\otimes\nu_1\oplus\nu_1)=(k+1)w_1(\ell)+2w_1(\nu_1)=0$. Noting that $q^*(\ell\otimes\nu_1)$ coincides with $q^*\nu_1$ this yields that we have
  \begin{alignat*}2
    2e(\ell\otimes\nu_1\oplus\nu_1)^{\frac r2-1}&=2e\left(\left(\frac r2-1\right)(\ell\otimes\nu_1\oplus\nu_1)\right)=q_!e\left(\left(\frac r2-1\right)(q^*(\ell\otimes\nu_1)\oplus q^*\nu_1)\right)=\\
                                                &=q_!e((r-2)q^*(\ell\otimes\nu_1))=2e((r-2)(\ell\otimes\nu_1))=2e(2(\ell\otimes\nu_1))^{\frac r2-1}=\\
                                                &=2p_{\frac{k+1}2}(\ell\otimes\nu_1)^{\frac r2-1}=2p_{\frac{k+1}2}(\nu_f\oplus\ell)^{\frac r2-1}=2p_{\frac{k+1}2}(\nu_f)^{\frac r2-1}.
  \end{alignat*}
  Thus the double of the oriented Thom polynomial evaluated on the map $f$ is
  $$2\Tp^\SO_{\Sigma^{1_r}(k)}(f)=(i_2)_!\left(2p_{\frac{k+1}2}(\nu_f)^{\frac r2-1}\right)=2(i_2)_!1_{\Sigma^{1,1}(f)}p_{\frac{k+1}2}(\nu_f)^{\frac r2-1}$$
  where $1_{\Sigma^{1,1}(f)}$ is the fundamental class of the manifold $\Sigma^{1,1}(f)$. But $2(i_2)_!1_{\Sigma^{1,1}(f)}$ is the double of the Thom polynomial $\Tp^\SO_{\Sigma^{1_2}(k)}$ evaluated on the map $f$ which, by \cite{singprim}, is $2p_{\frac{k+1}2}(\nu_f)$. This concludes the proof.
\end{prf}

\section{Open questions and concluding remarks}
\label{sec:fin}

The investigations in the present paper can be seen as working on the following (quite vague) problem:

\begin{ques}
  For which classes of smooth maps $f$ and which numbers $r$ do the cohomology classes $\Tp_{\Sigma^{1_r}(k)}(f)$ and $\Tp_{\Sigma^r(k-r+1)}(f)$ (or their oriented analogues) coincide? When is such a coincidence of Thom polynomials explained by a cobordism, bordism or other more geometric equivalence of the respective singular loci?
\end{ques}

For $r=1$ the (closures of the) singular loci in question coincide by definition for any $f$. For $r=2$ the coincidence also holds for any $f$ and the second question above in this case is the question of Fehér. We answered it in theorem \ref{thm:cusp} by an embedded cobordism in the case when $f$ is a cusp map. Although this already explains the coincidence of Thom polynomials, it is natural to ask the following:

\begin{ques}
  Does theorem \ref{thm:cusp} remain true if we replace the word ``cusp'' by ``Morin'' in it? Is there a similar equivalence (which is stronger than homology) of the $\Sigma^{1_2}$-locus of a generic map $f$ and the $\Sigma^2$-locus of its de-suspension?
\end{ques}



By a ``stronger equivalence'' here we mean any geometric equivalence of the singular loci noting that for a general map $f$ with de-suspension $\tilde f$ the closures $\ol\Sigma^{1_2}(f)$ and $\ol\Sigma^2(\tilde f)$ are not manifolds, hence we cannot talk about usual cobordism for them. Instead, we could ask for example whether they are cobordant as defective embeddings in some generalisation of the sense of Rost (see e.g. \cite{defekt}).

\begin{rmk}
  Such a stronger equivalence could also explain the fact that not only the Thom polynomials but even the avoiding ideals of the singularities $\Sigma^{1_2}(k)$ and $\Sigma^2(k-1)$ coincide, that is, the kernels of the restriction homomorphisms $H^*(B\O;\Z_2)\to H^*(K_{\{\Sigma^0(k),\Sigma^{1,0}(k)\}};\Z_2)$ and $H^*(B\O;\Z_2)\to H^*(K_{\{\Sigma^0(k-1),\Sigma^1(k-1)\}};\Z_2)$ which were computed in \cite{avcusp} and \cite{avsigmar} respectively.
\end{rmk}

On another line, theorem \ref{thm:prim} and corollary \ref{crly:prim} answer the first set of questions above in the case when $f$ is prim and theorem \ref{thm:twprim} and corollary \ref{crly:twprim} give a partial answer when $f$ is twisted prim. We conjecture that in theorem \ref{thm:twprim} and corollary \ref{crly:twprim} the factors $2$ can be omitted, more concretely, we pose the following:

\begin{conj}
  If $f\colon M^n\to N^{n+k}$ is a twisted prim map where $M$ is closed and $k\ge r-1$ and $\tilde f\colon M\times\R^{r-1}\to N$ is a de-suspension of $f$, then the submanifolds $\ol\Sigma^{1_r}(f)$ and $\Sigma^r(\tilde f)\cap M\times\{0\}$ of $M$ are cobordant as embedded submanifolds in $M$. This cobordism is oriented if $M$ and $N$ are oriented and the codimension $k$ is odd and $r$ is even.
\end{conj}

We can consider a stronger version of twisted prim maps: we say that a map $f\colon M\to N$ is \textit{strongly twisted prim} if there is a line bundle $\ell$ over $N$ (instead of $M$) such that $\ker df$ is $f^*\ell|_{\Sigma(f)}$. Then a similar consideration as the one implying that a Morin map with trivial kernel line bundle is prim, yields that there is an immersion $j$ from $M$ to the total space of $\ell$ which is a lift of $f$. For these strongly twisted prim maps, lemmas \ref{lemma:1r} and \ref{lemma:r} may extend using some modification of the normal bundle $\nu_j$ instead of $\nu_j$, thus giving a proof of the above conjecture in this case.

\begin{rmk}
  If the conjecture above is true, then we get obstructions for the existence of twisted prim maps: the Thom polynomials $\Tp_{\Sigma^{1_r}(k)}|_{\textrm{Morin}}$ and $\Tp^\SO_{\Sigma^{1_r}(k)}|_{\textrm{Morin}}$ 
  were computed in theorem \ref{thm:mor} while the Thom polynomials $\Tp_{\Sigma^r(k-r+1)}$ and $\Tp^\SO_{\Sigma^r(k-r+1)}$ are known from the Giambelli--Thom--Porteous formula and Ronga's computation \cite{calcsk}. The obstructions coming from this conjecture are the equalities $\Tp_{\Sigma^{1_r}(k)}(f)=\Tp_{\Sigma^r(k-r+1)}(f)$ and $\Tp^\SO_{\Sigma^{1_r}(k)}(f)=\Tp^\SO_{\Sigma^r(k-r+1)}(f)$ for twisted prim maps $f$.
\end{rmk}

\begin{rmk}
  A principle of Kazarian in \cite{morchar} says that when considering Thom polynomials of Morin maps, one can assume that they are, in our terminology, strongly twisted prim; that is, if a Thom polynomial formula holds for strongly twisted prim maps, then it is likely to hold for any Morin map. If Kazarian's principle also holds for the statement of the above conjecture, then we would get the same obstructions for the existence of general Morin maps which would uncover possible hidden relations in the avoiding ideal of $\Sigma^2(k)$. We note that an example for this principle is theorem \ref{thm:mor}: the formulae there can be proved for strongly twisted prim maps using the same ideas in a somewhat simpler way than for general Morin maps (in particular, there is no need to use lemma \ref{lemma:nu1}). But the actual proof we gave above shows that these ideas work for all Morin maps.
\end{rmk}

Lastly we note that lemmas \ref{lemma:1r}, \ref{lemma:r} and \ref{lemma:eq} also hold in the complex setting, hence the analogue of theorem \ref{thm:prim} and corollary \ref{crly:prim} for holomorphic prim maps is true as well. As we remarked in the introduction, the analogue of theorem \ref{thm:mor} is also true for holomorphic Morin maps, it was formulated and proved in \cite{morchar}. However, the proofs of theorems \ref{thm:cusp} and \ref{thm:twprim} do not go through for holomorphic maps (and the analogue of corollary \ref{crly:cusp} is not even true). Another line of investigation would be to try proving or disproving every statement in this paper for holomorphic maps instead of smooth maps, but, as it is usual with holomorphic problems, many smooth topological methods (e.g. relying on generic perturbations) would not work then.

\section*{Appendix: computational details in the proof of theorem \ref{thm:cusp}}
\label{sec:app}

In appropriate local coordinate systems around the points of $\Sigma^{1_2}(f)$ and their images, the map $\tilde f$ has the form $(p,t)\mapsto f(p)+t\sigma(p)$ where for simplicity we use the same letters for the maps written in local coordinates as for the original maps. Here $f\colon\R^n\to\R^{n+k}$ can be chosen as the normal form of the cusp singularity in these dimensions and $\sigma\colon\R^n\to\R^{n+k}$ maps a point $p\in\Sigma(f)$ to the vector $\frac{d^2f}{dq^2}(p)$ projected to the orthogonal complement of $\im df_p$ where $q$ generates the subspace $\ker df_p$.

The local form of cusp from $n$-dimensional space to $(n+k)$-dimensional is the $(n-2k-2)$-fold suspension of the prototype of the $k$-codimensional cusp which is, using \cite{mor}, given by
\begin{alignat*}2
  f\colon(x_1,\ldots,x_{2k},y,z,s_1,\ldots,s_{n-2k-2})&\mapsto(X_1,\ldots,X_{2k+1},Y_1,\ldots,Y_k,Z,S_1,\ldots,S_{n-2k-2})\\
  X_i&=x_i\quad(i\le2k),\\
  X_{2k+1}&=y,\\
  Y_i&=zx_{2i-1}+z^2x_{2i},\\
  Z&=zy+z^3,\\
  S_i&=s_i
\end{alignat*}
where $(x,y,z,s)=(x_1,\ldots,x_{2k},y,z,s_1,\ldots,s_{n-2k-2})$ and $(X,Y,Z,S)=(X_1,\ldots,X_{2k+1},Y_1,\allowbreak\ldots,Y_k,Z,S_1,\ldots,S_{n-2k-2})$ are the coordinates in $\R^n$ and $\R^{n+k}$ respectively and with a slight abuse of notation we also use the latter as the coordinate functions of $f$. Its derivative is
\bigskip
\begin{align*}
  df_{(x,y,z,s)}=\begin{NiceArray}{cccccc|c|c|ccc}[margin]
    \Block{6-1}{
      \scriptstyle X~~~} &1 \\
    && 1 \\
    &&& \ddots \\
    &&&& 1 \\
    &&&&& 1 \\
    &&&&&& 1 \\
    \cline{2-11}
    \Block{3-1}{
      \scriptstyle Y~~~} &z & z^2 &&&&& x_1+2zx_2 \\
    &&& \ddots &&&& \vdots \\
    &&&& z & z^2 && x_{2k-1}+2zx_{2k} \\
    \cline{2-11}
    \scriptstyle Z~~~ &&&&&& z & y+3z^2 \\
    \cline{2-11}
    \Block{3-1}{
      \scriptstyle S~~~} &&&&&&&& 1 \\
    &&&&&&&&& \ddots \\
    &&&&&&&&&& 1
    \CodeAfter
    \SubMatrix({1-2}{13-11})[left-xshift=2pt]
    \SubMatrix[{7-2}{7-3}][xshift=-2pt]
    \SubMatrix[{9-5}{9-6}][xshift=-2pt]
    \OverBrace[shorten,yshift=4pt]{1-2}{13-6}{
      \scriptstyle x}
    \OverBrace[shorten,yshift=4pt]{1-7}{1-7}{
      \scriptstyle y}
    \OverBrace[shorten,yshift=4pt]{1-8}{13-8}{
      \scriptstyle z}
    \OverBrace[shorten,yshift=4pt]{1-9}{13-11}{
      \scriptstyle s}
    \SubMatrix{\{}{1-2}{6-11}{.}[left-xshift=7pt]
    \SubMatrix{\{}{7-2}{9-11}{.}[left-xshift=7pt]
    \SubMatrix{\{}{10-2}{10-11}{.}[left-xshift=8pt]
    \SubMatrix{\{}{11-2}{13-11}{.}[left-xshift=7pt]
  \end{NiceArray}
\end{align*}
where empty spaces mean $0$ terms (this is a convention we will also use in the following). Hence $\Sigma(f)$ is given by the equations $x_{2i-1}=-2zx_{2i}$ and $y=-3z^2$ and here $\ker df$ is generated by the $z$-coordinate vector in $\R^n$ and $\im df$ is generated by the vectors $u_1,\ldots,u_{k+1},v_1,\ldots,v_k,w_1,\allowbreak\ldots,w_{n-2k-2}$ where we have
\begin{align*}
  u_i=
  \begin{NiceArray}{ccl}
    \Block{3-1}{\scriptstyle X~~} \\
    & 1 & ~\scriptstyle X_{2i-1}\\
    \\
    \Block{3-1}{\scriptstyle Y~~} \\
    & z & ~\scriptstyle Y_i\\
    \\
    \scriptstyle Z~~\\
    \Block{2-1}{\scriptstyle S~~}\\
    \\
    \CodeAfter
    \SubMatrix({1-2}{9-2})[hlines={3,6,7}]
    \SubMatrix{\{}{1-2}{3-2}{.}[left-xshift=5pt,extra-height=-2pt]
    \SubMatrix{\{}{4-2}{6-2}{.}[left-xshift=5pt,extra-height=-2pt]
    \SubMatrix{\{}{7-2}{7-2}{.}[left-xshift=6pt,extra-height=-2pt]
    \SubMatrix{\{}{8-2}{9-2}{.}[left-xshift=5pt,extra-height=-2pt]
    \SubMatrix{.}{2-2}{2-2}{\}}[right-xshift=6pt]
    \SubMatrix{.}{5-2}{5-2}{\}}[right-xshift=6pt]
  \end{NiceArray} (i\le k),
  \quad u_{k+1}=
  \begin{NiceArray}{ccl}
    \Block{3-1}{\scriptstyle X~~} \\
    \\
    & 1 & ~\scriptstyle X_{2k+1}\\
    \Block{2-1}{\scriptstyle Y~~} \\
    \\
    \scriptstyle Z~~ & z\\
    \Block{2-1}{\scriptstyle S~~}\\
    \\
    \CodeAfter
    \SubMatrix({1-2}{8-2})[hlines={3,5,6}]
    \SubMatrix{\{}{1-2}{3-2}{.}[left-xshift=5pt,extra-height=-2pt]
    \SubMatrix{\{}{4-2}{5-2}{.}[left-xshift=5pt,extra-height=-2pt]
    \SubMatrix{\{}{6-2}{6-2}{.}[left-xshift=6pt,extra-height=-2pt]
    \SubMatrix{\{}{7-2}{8-2}{.}[left-xshift=5pt,extra-height=-2pt]
    \SubMatrix{.}{3-2}{3-2}{\}}[right-xshift=6pt]
  \end{NiceArray}
  \text{and}\quad v_i=
  \begin{NiceArray}{ccl}
    \Block{3-1}{\scriptstyle X~~} \\
    & 1 & ~\scriptstyle X_{2i}\\
    \\
    \Block{3-1}{\scriptstyle Y~~} \\
    & z^2 & ~\scriptstyle Y_i\\
    \\
    \scriptstyle Z~~\\
    \Block{2-1}{\scriptstyle S~~}\\
    \\
    \CodeAfter
    \SubMatrix({1-2}{9-2})[hlines={3,6,7}]
    \SubMatrix{\{}{1-2}{3-2}{.}[left-xshift=5pt,extra-height=-2pt]
    \SubMatrix{\{}{4-2}{6-2}{.}[left-xshift=5pt,extra-height=-2pt]
    \SubMatrix{\{}{7-2}{7-2}{.}[left-xshift=6pt,extra-height=-2pt]
    \SubMatrix{\{}{8-2}{9-2}{.}[left-xshift=5pt,extra-height=-2pt]
    \SubMatrix{.}{2-2}{2-2}{\}}[right-xshift=6pt]
    \SubMatrix{.}{5-2}{5-2}{\}}[right-xshift=6pt]
  \end{NiceArray}
\end{align*}
and $w_i$ is the $S_i$-coordinate vector in $\R^{n+k}$. Now the second intrinsic derivative defining the section $\sigma$ is the projection of the vector
\begin{align*}
  \frac{d^2f}{dz^2}=
  \begin{NiceArray}{cc}
    \Block{2-1}{\scriptstyle X~~} \\
    \\
    \Block{3-1}{\scriptstyle Y~~} & 2x_2\\
    & \vdots \\
    & 2x_{2k} \\
    \scriptstyle Z~~ & 6z\\
    \Block{2-1}{\scriptstyle S~~}\\
    \\
    \CodeAfter
    \SubMatrix({1-2}{8-2})[hlines={2,5,6}]
    \SubMatrix{\{}{1-2}{2-2}{.}[left-xshift=5pt,extra-height=-2pt]
    \SubMatrix{\{}{3-2}{5-2}{.}[left-xshift=5pt,extra-height=-2pt]
    \SubMatrix{\{}{6-2}{6-2}{.}[left-xshift=6pt,extra-height=-2pt]
    \SubMatrix{\{}{7-2}{8-2}{.}[left-xshift=5pt,extra-height=-2pt]
  \end{NiceArray}
\end{align*}
to the orthogonal complement of the subspace generated by $u_1,\ldots,u_{k+1},v_1,\ldots,v_k$ (since it is already orthogonal to the vectors $w_i$). Although this only defines $\sigma$ restricted to $\Sigma(f)$, we may assume that the same formula gives $\sigma$ on a neighbourhood of it as well. To compute this formula, we first project $v_i$ to the orthogonal complement of $u_i$ for all $i$ which yields
\begin{align*}
  \tilde v_i=v_i-\frac{\la u_i,v_i\ra}{||u_i||^2}u_i=
  \begin{NiceArray}{ccl}
    \Block{4-1}{\scriptstyle X~~} \\
    & -\frac{z^3}{1+z^2} & ~\scriptstyle X_{2i-1} \\
    & 1 & ~\scriptstyle X_{2i} \\
    \\
    \Block{3-1}{\scriptstyle Y~~} \\
    & \frac{z^2}{1+z^2} & ~\scriptstyle Y_i \\
    \\
    \scriptstyle Z~~ \\
    \Block{2-1}{\scriptstyle S~~}\\
    \\
    \CodeAfter
    \SubMatrix({1-2}{10-2})[hlines={4,7,8}]
    \SubMatrix{\{}{1-2}{4-2}{.}[left-xshift=5pt,extra-height=-2pt]
    \SubMatrix{\{}{5-2}{7-2}{.}[left-xshift=5pt,extra-height=-2pt]
    \SubMatrix{\{}{8-2}{8-2}{.}[left-xshift=6pt,extra-height=-2pt]
    \SubMatrix{\{}{9-2}{10-2}{.}[left-xshift=5pt,extra-height=-2pt]
    \SubMatrix{.}{2-2}{2-2}{\}}[right-xshift=6pt,extra-height=-2pt]
    \SubMatrix{.}{3-2}{3-2}{\}}[right-xshift=6pt,extra-height=-2pt]
    \SubMatrix{.}{6-2}{6-2}{\}}[right-xshift=6pt,extra-height=-2pt]
  \end{NiceArray}
\end{align*}
and now $u_1,\ldots,u_{k+1},\tilde v_1,\ldots,\tilde v_k$ is a pairwise orthogonal basis for the subspace generated by $u_1,\ldots,u_{k+1},v_1,\ldots,v_k$. Hence we have
\begin{align*}
  \sigma=\frac{d^2f}{dz^2}-\sum_{i=1}^{k+1}\frac{\la u_i,\frac{d^2f}{dz^2}\ra}{||u_i||^2}u_i-\sum_{i=1}^k\frac{\la \tilde v_i,\frac{d^2f}{dz^2}\ra}{||\tilde v_i||^2}\tilde v_i
\end{align*}
whose $X_{2i-1},X_{2i},X_{2k+1},Y_i$ and $Z$ coordinates (for $i\le k$) are $-\frac{2zx_{2i}}{1+z^2+z^4},-\frac{2z^2x_{2i}}{1+z^2},-\frac{6z^3}{1+z^2},\allowbreak\frac{2x_{2i}}{1+z^2+z^4}$ and $\frac{6z}{1+z^2}$ respectively. Thus the map $\tilde f\colon\R^{n+1}\to\R^{n+k}$ defined by $\tilde f(p,t):=f(p)+t\sigma(p)$ is of the form
\begin{alignat*}2
  \tilde f\colon(x_1,\ldots,x_{2k},y,z,s_1,\ldots,s_{n-2k-2},t)&\mapsto(\tilde X_1,\ldots,\tilde X_{2k+1},\tilde Y_1,\ldots,\tilde Y_k,\tilde Z,\tilde S_1,\ldots,\tilde S_{n-2k-2})\\
  \tilde X_{2i-1}&=x_{2i-1}-\frac{2tzx_{2i}}{1+z^2+z^4}\quad(i\le2k),\\
  \tilde X_{2i}&=x_{2i}-\frac{2tz^2x_{2i}}{1+z^2},\\
  \tilde X_{2k+1}&=y-\frac{6tz^3}{1+z^2},\\
  \tilde Y_i&=zx_{2i-1}+z^2x_{2i}+\frac{2tx_{2i}}{1+z^2+z^4},\\
  \tilde Z&=zy+z^3+\frac{6tz}{1+z^2},\\
  \tilde S_i&=s_i
\end{alignat*}
and so its derivative $d\tilde f$ at $(x,y,z,s,t)$ is
\bigskip
\begin{align*}
  \begin{NiceArray}{cccccc|c|c|ccc|c}[cell-space-limits=2pt,margin]
    \Block{6-1}{\scriptstyle X~~~} & 1 & -\frac{2tz}{1+z^2+z^4} &&&&& \frac{2tx_2(z^2+3z^4-1)}{(1+z^2+z^4)^2} &&&& -\frac{2zx_2}{1+z^2+z^4} \\
    & 0 & 1-\frac{2tz^2}{1+z^2} &&&&& -\frac{4tzx_2}{(1+z^2)^2} &&&& -\frac{2z^2x_2}{1+z^2} \\
    &&& \ddots &&&& \vdots &&&& \vdots \\
    &&&& 1 & -\frac{2tz}{1+z^2+z^4} && \frac{2tx_{2k}(z^2+3z^4-1)}{(1+z^2+z^4)^2} &&&& -\frac{2zx_{2k}}{1+z^2+z^4} \\
    &&&& 0 & 1-\frac{2tz^2}{1+z^2} && -\frac{4tzx_{2k}}{(1+z^2)^2} &&&& -\frac{2z^2x_{2k}}{1+z^2} \\
    &&&&&& 1 & -\frac{6tz^2(z^2+3)}{(1+z^2)^2} &&&& -\frac{6z^3}{1+z^2} \\
    \cline{2-12}
    \Block{3-1}{\scriptstyle Y~~~} & z & z^2+\frac{2t}{1+z^2+z^4} &&&&& x_1+2zx_2\left(1-\frac{t(1+z^2)}{(1+z^2+z^4)^2}\right) &&&& \frac{2x_2}{1+z^2+z^4} \\
    &&& \ddots &&&& \vdots &&&& \vdots \\
    &&&& z & z^2+\frac{2t}{1+z^2+z^4} && x_{2k-1}+2zx_{2k}\left(1-\frac{t(1+z^2)}{(1+z^2+z^4)^2}\right) &&&& \frac{2x_{2k}}{1+z^2+z^4} \\
    \cline{2-12}
    \scriptstyle Z~~~ &&&&&& z & y+3z^2+\frac{6t(1-z^3)}{(1+z^2)^2} &&&& \frac{6z}{1+z^2} \\
    \cline{2-12}
    \Block{3-1}{\scriptstyle S~~~} &&&&&&&& 1 \\
    &&&&&&&&& \ddots \\
    &&&&&&&&&& 1
    \CodeAfter
    \SubMatrix({1-2}{13-12})[left-xshift=2pt]
    \SubMatrix[{1-2}{2-3}][right-xshift=-3pt,left-xshift=-1pt]
    \SubMatrix[{4-5}{5-6}][right-xshift=-3pt,left-xshift=-1pt]
    \SubMatrix[{1-8}{2-8}][xshift=-35pt]
    \SubMatrix[{4-8}{5-8}][xshift=-35pt]
    \SubMatrix[{1-12}{2-12}][xshift=-2pt]
    \SubMatrix[{4-12}{5-12}][xshift=-2pt]
    \SubMatrix[{7-2}{7-3}][xshift=-1pt,extra-height=-3pt]
    \SubMatrix[{9-5}{9-6}][xshift=-1pt,extra-height=-3pt]
    \OverBrace[shorten,yshift=4pt]{1-2}{13-6}{\scriptstyle x}
    \OverBrace[shorten,yshift=4pt]{1-7}{1-7}{\scriptstyle y}
    \OverBrace[shorten,yshift=4pt]{1-8}{13-8}{\scriptstyle z}
    \OverBrace[shorten,yshift=4pt]{1-9}{13-11}{\scriptstyle s}
    \OverBrace[shorten,yshift=4pt]{1-12}{1-12}{\scriptstyle t}
    \SubMatrix{\{}{1-2}{6-12}{.}[left-xshift=7pt]
    \SubMatrix{\{}{7-2}{9-12}{.}[left-xshift=7pt]
    \SubMatrix{\{}{10-2}{10-12}{.}[left-xshift=7pt,extra-height=-3pt]
    \SubMatrix{\{}{11-2}{13-12}{.}[left-xshift=7pt]
  \end{NiceArray}~.
\end{align*}

Now we need that the map $d\tilde f\colon\R^{n+1}\to\Hom(\R^{n+1},\R^{n+k})$ is transverse to the subspace $\Sigma^2(\R^{n+1},\R^{n+k})$ consisting of corank-$2$ homomorphisms. The normal space of $\Sigma^2(\R^{n+1},\R^{n+k})$ in $\Hom(\R^{n+1},\R^{n+k})$ at a point $\varphi\colon\R^{n+1}\to\R^{n+k}$ is $\Hom(\ker\varphi,\coker\varphi)$; we want to show that the projection of the image of $d(d\tilde f)_{p}$ to this subspace is surjective for all points $p\in(d\tilde f)^{-1}(\Sigma^2(\R^{n+1},\R^{n+k}))=\Sigma^2(\tilde f)$.

The singular locus $\Sigma^2(\tilde f)$ is the subspace $(x,y,z)=(0,0,0)$ in $\R^{n+1}$ and for any point $\varphi\in\im d\tilde f|_{\Sigma^2(\tilde f)}$ the subspace $\Hom(\ker\varphi,\coker\varphi)<\Hom(\R^{n+1},\R^{n+k})$ is generated by the coordinate vectors corresponding to $yY_1,\ldots,yY_k,yZ,\allowbreak tY_1,\ldots,tY_k,tZ$ using the notation that for each coordinate $q$ of $\R^{n+1}$ and $Q$ of $\R^{n+k}$ we denote by $qQ$ the corresponding coordinate of $\Hom(\R^{n+1},\R^{n+k})$. It is not hard to see that every such coordinate vector is indeed in the image of $d(d\tilde f)_{(0,0,0,s,t)}$ projected to it (for any $s,t$), namely in the partial derivatives $\frac{d(d\tilde f)}{dx_{2i-1}}(0,0,0,s,t)$, $\frac{d(d\tilde f)}{dy}(0,0,0,s,t)$, $\frac{d(d\tilde f)}{dx_{2i}}(0,0,0,s,t)$ and $\frac{d(d\tilde f)}{dz}(0,0,0,s,t)$ respectively the coordinate $yY_i$, $yZ$, $tY_i$ and $tZ$ is non-zero and all other coordinates of this subspace are $0$. This is what we wanted to show.


\begin{thebibliography}{0000000}{\footnotesize
    
  \bibitem[BSz12]{tpmor} G. Bérczi, A. Szenes, \textit{Thom polynomials of Morin singularities}, Ann. of Math. 175 (2012), 567--629.

  \bibitem[Bo67]{singdiff}  J. M. Boardman, {\it Singularities of differentiable maps}, Publ. Math. I.H.É.S. 33 (1967), 21--57.
    
  \bibitem[BH61]{borhae} A. Borel, A. Haefliger, {\it La classe d'homologie fondamentale d'un espace analytique}, Bull. Soc. Math. France 89 (1961), 461--513.

  \bibitem[Da74]{damonflag} J. Damon, \textit{The Gysin homomorphism for flag bundles: applications}, Amer. J. Math. 96 (1974), 248--260.
    
  \bibitem[FR02]{inttp} L. M. Fehér, R. Rimányi, \textit{Thom polynomials with integer coefficients}, Illinois J. Math. 46 (2002), 1145--1158.
    
  \bibitem[FR04]{calctp} L. M. Fehér, R. Rimányi, \textit{Calculation of Thom polynomials and other cohomological obstructions for group actions}, in: Real and Complex Singularities (São Carlos, 2002), Contemp. Math. 354, Amer. Math. Soc. (2004), 69--93.
    
  \bibitem[Ka97]{charsing} M. É. Kazarian, \textit{Characteristic classes of singularity theory}, in: The Arnold--Gelfand Mathematical Seminars, Birkhäuser (1997), 325--340.

  \bibitem[Ka01]{kazspace} M. É. Kazarian, {\it Classifying spaces of singularities and Thom polynomials}, in: New Developments in Singularity Theory (Cambridge, 2000), NAII 21, Springer (2001), 117--134.

  \bibitem[Ka06]{kaztp} M. É. Kazarian, {\it Thom polynomials}, Singularity theory and its applications, Adv. Stud. Pure Math. 43 (2006), 85--135.
    
  \bibitem[Ka]{morchar} M. É. Kazarian, \textit{Morin maps and their characteristic classes}, unpublished preprint.


    
  \bibitem[MS74]{charclass} J. W. Milnor, J. D. Stasheff, \textit{Characteristic classes}, Ann. of Math. Stud. 76, Princeton University Press (1974).
    
  \bibitem[Mo65]{mor} B. Morin, {\it Formes canoniques des singularités d'une application différentiable}, C. R. Acad. Sci. Paris 260 (1965), 5662--5665 and 6503--6506.
    
  \bibitem[Po71]{gtp} I. R. Porteous, \textit{Simple singularities of maps}, in: Proc. Liverpool Singularities I, Lecture Notes in Math. 192, Springer (1971), 286--307.
    
  \bibitem[Pr88]{avsigmar} P. Pragacz, \textit{Enumerative geometry of degeneracy loci}, Ann. Sci. École Norm. Sup. (4) 21 (1988) 413--454.
    
  \bibitem[Ri00]{orsing} R. Rimányi, \textit{On the orientability of singularity submanifolds}, J. Math. Soc. Japan 52 (2000), 91--98.
    
  \bibitem[Ri01]{rim} R. Rimányi, {\it Thom polynomials, symmetries and incidences of singularities}, Invent. Math. 143 (2001), 499--521.
    
  \bibitem[Ro71]{calcsk} F. Ronga, \textit{Le calcul de la classe de cohomology entière duale à $\ol\Sigma^k$}, in: Proc. Liverpool Singularities I, Lecture Notes in Math. 192, Springer (1971), 313--315.
    
  \bibitem[Ro72]{calcord2} F. Ronga, \textit{Le calcul des classes duales aux singularités de Boardman d'ordre deux}, Comment. Math. Helv. 47 (1972), 15--35.
    
  \bibitem[Ro]{tps11} F. Ronga, \textit{The integral Thom polynomial for $\Sigma^{1,1}$}, unpublished preprint.

  \bibitem[Ros87]{defekt} M. Rost, \textit{Abbildungsdefekte in $4$-Mannigfaltigkeiten}, Regensburger Math. Schriften 12, Universität Regensburg (1987).

  \bibitem[Sz91]{immemb} A. Szűcs, \textit{On the cobordism groups of immersions and embeddings}, Math. Proc. Camb. Phil. Soc. 109 (1991), 343--349.
    
  \bibitem[Sz00]{singprim} A. Szűcs, \textit{On the singularities of hyperplane projections of immersions}, Bull. London Math. Soc. 32 (2000), 364--374.

  \bibitem[Sz08]{hosszu} A. Szűcs, {\it Cobordism of singular maps}, Geom. Topol. 12 (2008), 2379--2452.
    
  \bibitem[Te09]{avcusp} T. Terpai, \textit{Calculation of the avoiding ideal for $\Sigma^{1,1}$}, in: Algebraic topology---old and new, Banach Center Publ. 85 (2009), 307--313.
    
  \bibitem[Th56]{tp} R. Thom, {\it Les singularités des applications différentiables}, Ann. Inst. Fourier 6 (1956), 43--87.
    
  \bibitem[Va88]{vas} V. A. Vassiliev, \textit{Lagrange and Legendre characteristic classes}, Adv. Stud. Contemp. Math. 3, Gordon and Breach (1988).
    
}\end{thebibliography}
\end{document}